\documentclass[twoside]{IEEEtran}

\usepackage{graphicx}
\usepackage{booktabs}
\usepackage{multirow}
\usepackage{dsfont}

\usepackage{color}
\usepackage{cite}
\usepackage{amsfonts,amsmath,amssymb,amsthm}
\usepackage[utf8]{inputenc}
\usepackage[normalem]{ulem}

\newtheorem{definition}{\textbf{Definition}}

\newtheorem{theorem}{\textbf{Theorem}}

\DeclareMathOperator*{\esssup}{ess\,sup}
\newcommand{\Prob}{\mathbb{P}}
\newcommand{\Expect}{\mathbb{E}}
\newcommand{\ADD}{{\mathsf{ADD}}}
\newcommand{\CADD}{{\mathsf{CADD}}}
\newcommand{\WADD}{{\mathsf{WADD}}}
\newcommand{\FAR}{{\mathsf{FAR}}}
\newcommand{\FDR}{{\mathsf{FDR}}}
\newcommand{\ARL}{{\mathsf{ARL}}}
\newcommand{\PFA}{\mathsf{PFA}}
\newcommand{\mc}[1]{\mathcal{#1}}
\newcommand{\Dd}{{\mc{D}}}

\newcommand{\taus}{\tau_{\scriptscriptstyle \text{S}}}
\newcommand{\taush}{\tau_{\scriptscriptstyle \text{Sh}}}
\newcommand{\taum}{\tau_{\scriptscriptstyle \text{M}}}
\newcommand{\tauc}{\tau_{\scriptscriptstyle \text{C}}}
\newcommand{\tausr}{\tau_{\scriptscriptstyle \text{SR}}}

\newcommand{\taug}{\tau_{\scriptscriptstyle \text{G}}}

\title{Sequential (Quickest) Change Detection: \\Classical Results and New Directions}

\author{Liyan Xie,~\IEEEmembership{Student Member,~IEEE}, Shaofeng Zou,~\IEEEmembership{Member,~IEEE},\\ Yao Xie,~\IEEEmembership{Member,~IEEE}, and Venugopal V. Veeravalli,~\IEEEmembership{Fellow,~IEEE} \thanks{Manuscript received October 28, 2020; revised January 28, 2021; accepted April 5, 2021. V. V. Veeravalli was supported in part by the Army Research Laboratory under Cooperative Agreement W911NF-17-2-0196 (IoBT
CRA), and in part by the National Science Foundation (NSF) under grant CCF 16-18658 and CIF 15-14245, through the University of Illinois at Urbana-Champaign. L. Xie and Y. Xie were partially supported by an NSF CAREER Award CCF-1650913, DMS-1938106, DMS-1830210, CCF-1442635, and CMMI-1917624. S. Zou was partially supported by an NSF Award CCF-1948165.

L. Xie and Y. Xie are with H. Milton Stewart School of Industrial and Systems Engineering, Georgia Institute of Technology, Atlanta, GA 30332 USA (email: lxie49@gatech.edu; yao.xie@isye.gatech.edu). 

S. Zou is with the Department of Electrical Engineering, University at Buffalo, The State University of New York, Buffalo, NY 14260 USA (email: szou3@buffalo.edu). 

V. V. Veeravalli is with the Department of Electrical and Computer Engineering, University of Illinois at Urbana-Champaign, Urbana, IL 61801 USA (email: vvv@illinois.edu).}}

\begin{document}

\maketitle

\begin{abstract}
Online detection of changes in stochastic systems, referred to as sequential change detection or quickest change detection, is an important research topic in statistics, signal processing, and information theory, and has a wide range of applications. This survey starts with the basics of sequential change detection, and then moves on to generalizations and extensions of sequential change detection theory and methods. We also discuss some new dimensions that emerge at the intersection of sequential change detection with other areas, along with a selection of modern applications and remarks on open questions.
\end{abstract}


\section{Introduction}

\IEEEPARstart{T}{he} efficient detection of abrupt changes in the statistical behavior of streaming data is a classical and fundamental problem in signal processing and statistics. The abrupt change-point usually corresponds to a triggering event that could have catastrophic consequences if it is not detected in a timely manner. Therefore, the goal is to detect the change as quickly as possible, subject to false alarm constraints. Such problems have been studied under the theoretical framework of sequential (or quickest) change detection \cite{poor-hadj-QCD-book-2008,tartakovsky2014sequential,veeravalli2013quickest}.
With an increasing availability of high-dimensional streaming data, sequential change  detection has become a centerpiece for many real-world applications, ranging from monitoring power networks \cite{chen2015quickest}, internet traffic \cite{LakhinaCrovellaDiot2004},  cyber-physical systems \cite{cyberphysical15}, sensor networks \cite{VenuSensor2010}, social networks \cite{raginsky_OCP,PeelClauset2014}, epidemic detection \cite{epidemic_change_04}, scientific imaging \cite{solarFlares}, genomic signal processing \cite{ShenZhang12}, seismology \cite{change-point-seismology07}, video surveillance \cite{LeeKriegman2005}, and wireless communications \cite{LaiCognitive08}.  

In various applications, the streaming data is high-dimensional and collected over networks, such as social networks,  sensor networks, and cyber-physical systems. For this reason, the modern sequential change detection problem's scope has been extended far beyond its traditional setting, often challenging the assumptions made by classical methods. These challenges include complex spatial and temporal dependence of the data streams, transient and dynamic changes, high-dimensionality, and structured changes, as explained below. These challenges have fostered new advances in sequential change detection theory and methods in recent years.

(1) {\it Complex data distributions.} 
In modern applications, sequential data could have a complex spatial and temporal dependency, for instance, induced by the network structure  \cite{cyber_physical10,NetworkLasso2015,change-point-correlation-networks}. In social networks, dependencies are usually due to interaction and information diffusion \cite{li2017detecting}: users in the social network have behavior patterns influenced by their past, while at the same time, each user in the network will be influenced by friends and connections. In sensor networks for river contamination monitoring \cite{chen2017textsf}, sensor observations tend to be spatially and temporally correlated. 

\vspace{0.05in}
(2) {\it Data dynamics.} The statistical behavior of sequential data is often non-stationary, particularly in the post-change regime due to the dynamic behavior of the anomaly that causes the change. For example, after a linear outage in the power systems, the system's transient behavior is dominated by the generators' inertial response, and the post-change statistical behavior can be modeled using a sequence of temporally cascaded transient phases \cite{Rovatsos2:2016}.

\vspace{0.05in}
(3) {\it High-dimensionality.} Sequential data in modern applications is usually high-dimensional. For example, in sensor networks, the Long Beach 3D seismic array consists of approximately 5300 seismic sensors that record data continuously for seismic activity detection and analysis. Changes in high-dimensional time series usually exhibit low-dimensional structures in the form of sparsity, low-rankness, and subset structures, which can be exploited to enhance the capability to detect weak signals quickly. 

In this tutorial, our aim is to introduce standard methods and fundamental results in sequential change detection, along with recent advances. We also present new dimensions at the intersection of sequential change detection with other areas, as well as a selection of modern applications. We should emphasize that our focus is on {\it sequential} change detection, where the goal is to detect the change from sequential data in {\it real-time} and {\it as soon as possible}. Another important line of related research is {\it offline} change detection (e.g., \cite{siegmund2011detecting,frick2014multiscale}), where the goal is to identify and localize changes in data sequence in a {\it  retrospective} manner, which is not our focus here. Prior books and surveys on related topics include, for instance, change detection for dynamic systems \cite{lai1995sequential}, sequential analysis \cite{lai2001sequential,tartakovsky2014sequential}, sequential change detection \cite{veeravalli2013quickest,poor-hadj-QCD-book-2008,detectAbruptChange93}, Bayesian change detection \cite{tartakovsky2010state}, change detection assuming known pre- and post-change distributions \cite{polunchenko2012state} and using likelihood-based approaches \cite{siegmund2013change}, as well as time-series change detection \cite{aminikhanghahi2017survey}.

The rest of the survey is organized as follows. In Section \ref{sec:classic}, we present the basic problem setup and classical results. In Section \ref{sec:new-advances}, we discuss several extensions and generalizations of the classical methods. In Section \ref{sec:new-dimension}, we discuss new dimensions which intersect with sequential change detection, with some remarks on open questions. In Section \ref{sec:applications}, we present some modern applications of sequential change detection. In Section \ref{sec:discussion}, we make some concluding remarks.

\section{Classical Results} \label{sec:classic}

\subsection{Problem Definition}\label{sec:pro_def}

In the sequential change detection problem, also known as the quickest change detection (QCD) problem \cite{moulin2018statistical,poor-hadj-QCD-book-2008,veeravalli2013quickest}, the aim is to detect a possible change in the data generating distribution of a sequence of observations  $\{X_n, n=1, 2, \ldots\}$. The initial distribution of the observations is the one corresponding to normal system operation. At some unknown time $\gamma$ (referred to as the {\it change-point}), due to some event, the distribution of the random observations changes.  The goal is to detect the change as quickly as possible, subject to false-alarm constraints. We start by assuming that the observations are independent and identically distributed (i.i.d.) with probability density function (pdf) $f_0$ before and pdf $f_1$ after the change-point, respectively. We discuss generalizations to non-i.i.d. observations in Section \ref{sec:new-advances}.

\begin{figure}[!ht]
\centering
\begin{tabular}{cc}
\includegraphics[width=0.45\linewidth]{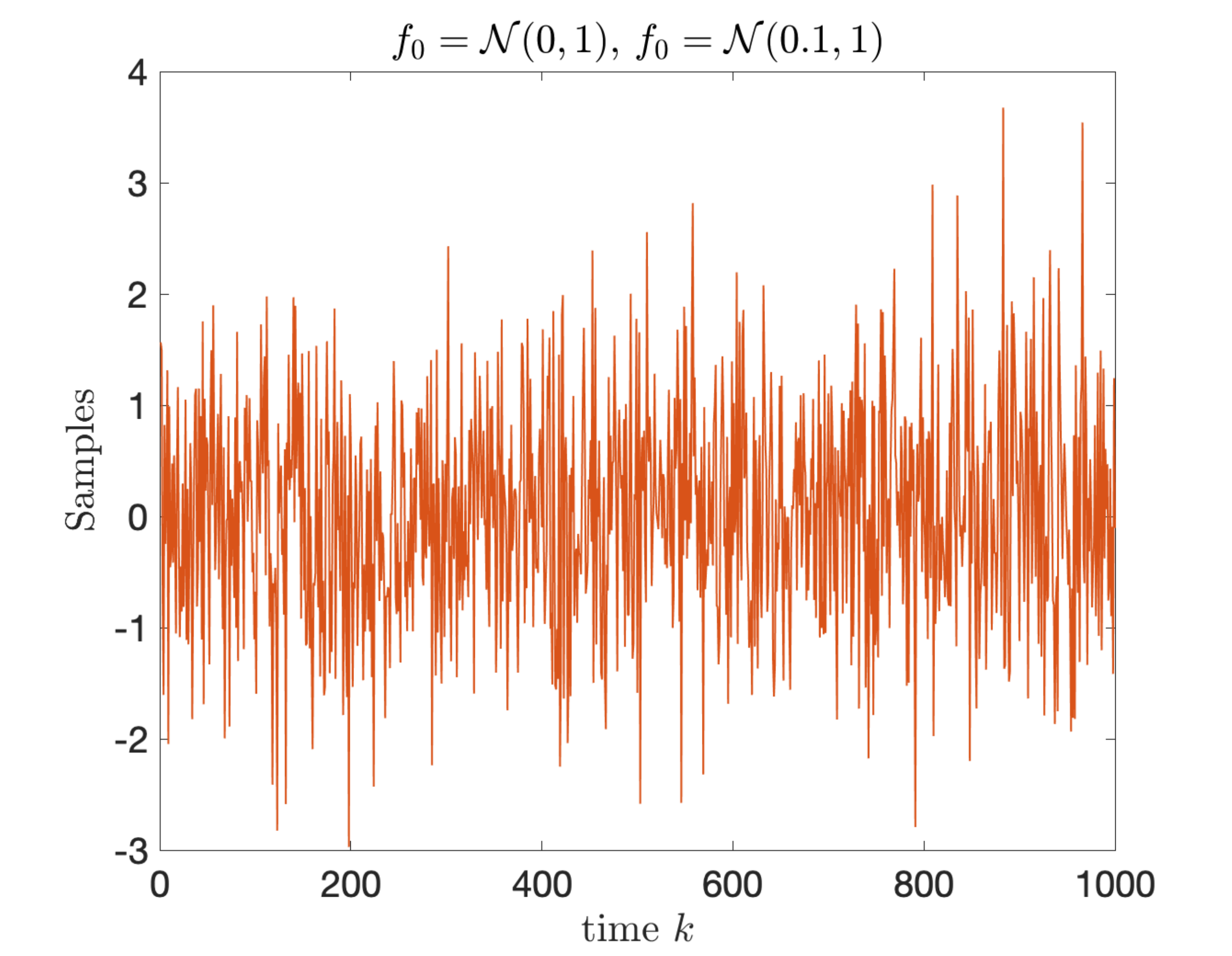} &
\includegraphics[width=0.45\linewidth]{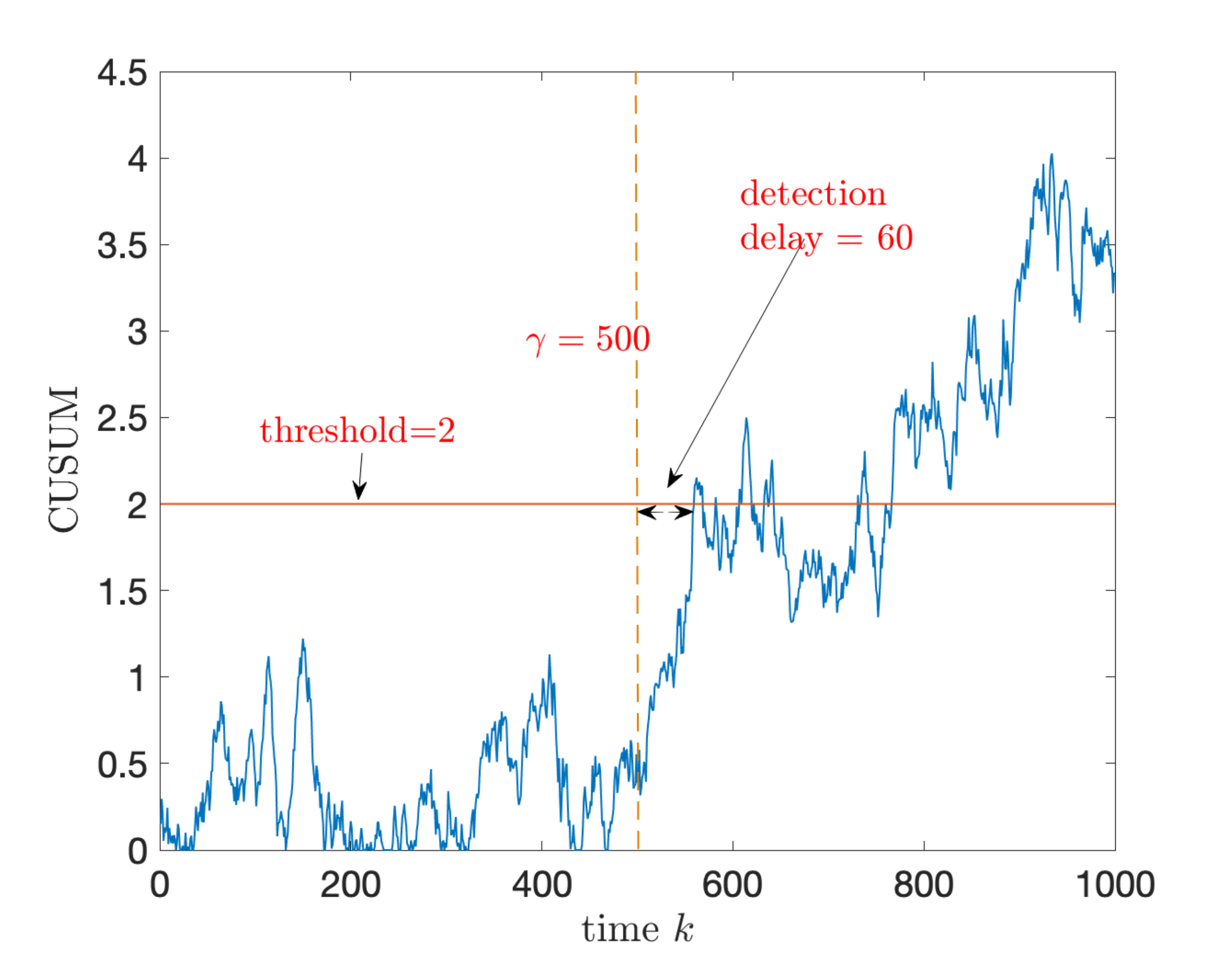}\\
(a) & (b)
\end{tabular}
\caption{To motivate the need for sequential change detection procedures, we plot a sample path with samples distributed according to $\mathcal N(0, 1)$ before the change and $\mathcal N(0.1, 1)$ after the change. We set the change-point $\gamma = 500$. As illustrated in (a), such a small mean shift cannot be detected through manual inspection of the samples. In (b), we plot the evolution of the CUSUM algorithm (detailed in Section~\ref{sec:CUSUM}) corresponding to the observations in (a), which can detect the change quickly.}
\label{fig:change-point-example}
\end{figure}

To motivate the design of algorithms for sequential change detection, we consider the example of detecting a change in the mean of the data generating distribution. In Fig.\,\ref{fig:change-point-example}(a), we plot a sample path of observations that are distributed according to a normal distribution with zero mean and unit variance $\mathcal N(0, 1)$ before the change-point of 500, and $\mathcal N(0.1, 1)$ after the change-point. As can be seen in Fig.\,\ref{fig:change-point-example}(a), such a small mean shift cannot be detected through manual inspection of the samples. In Fig.\,\ref{fig:change-point-example}(b), we plot the evolution path of a sequential change detection procedure, the CUSUM algorithm (which is discussed in detail in Section \ref{sec:CUSUM}), applied to the observations in Fig.\,\ref{fig:change-point-example}(a). It can be seen that the test statistic stays close to zero before the change and has a positive drift after the change. Therefore, the change can be detected by comparing the test statistic to a positive threshold $b$ (for instance, $b=2$) and raising an alarm when the test statistic exceeds the threshold for the first time. For the sample path in Fig.\,\ref{fig:change-point-example}(a), this approach incurs a detection delay of 60 samples (if we take samples daily, this means a detection delay of 2 months; if the sampling rate is 60 samples per second, this means a detection delay of one second). One natural question to ask is that: can we do better, at least on average? Clearly, if we set a lower threshold, for instance $b=1$, we can detect the change much more quickly. However, this would result in a false alarm at $k=112$. This example illustrates the {\it tradeoff} between false-alarm and detection delay, which is a central problem when designing sequential change detection procedures. The goal in sequential change detection theory is to find detection procedures that have guaranteed optimality properties in terms of this tradeoff.

\subsection{Mathematical Preliminaries}
\label{sec:preliminaries}

Sequential change detection is closely related to the problem of statistical hypothesis testing, in which observations, whose distribution depends on the hypothesis, are used to decide which of the hypotheses is true. For the special case of binary hypothesis testing, we have two hypotheses, the {\it null} hypothesis and the {\it alternative} hypothesis. The classic Neyman-Pearson Lemma \cite{neyman1933ix} establishes the form of the optimal test for this problem. In particular, consider the case of a single observation $X$, and suppose the pdf of $X$ under the null and alternative hypotheses are $f_0$ and $f_1$,  respectively. Then, the test that minimizes the false negative error (Type-II error), under the constraint of the false positive error (Type-I error), is to compare the {\em likelihood ratio} $f_1(X)/f_0(X)$ to a threshold to decide which hypothesis is true. The likelihood ratio test is also optimal under other criteria such as Bayesian and minimax \cite{moulin2018statistical}. As we will see, the likelihood ratio also plays a key role in the development of sequential change detection algorithms.

The goal of sequential change detection is to design a {\em stopping time} on the observation sequence at which it is declared that a change has occurred. A stopping time is formally defined as follows:
\begin{definition}[Stopping Time] A  stopping time with respect to a random sequence $\{X_n, n=1, 2, \ldots\}$ is a random
variable $\tau$ such that for each $n$, the event $\{ \tau=n\} \in \sigma (X_1, \ldots, X_n)$,
where $\sigma (X_1, \ldots, X_n)$ denotes the sigma-algebra generated by $(X_1, \ldots, X_n)$.
Equivalently, the event $\{ \tau =n\}$  is a function of only $X_1, \ldots, X_n$.
\end{definition}
The main results on stopping times that are most useful for sequential change detection problems include Doob's Optional Stopping Theorem \cite{chow-robb-sieg-book-1971} and Wald's Identity \cite{Siegmund1985}.  

A quantity that plays an important role in the performance of sequential change detection algorithms is the Kullback-Leibler (KL) divergence between two distributions.
\begin{definition} \label{def:KL-Divergence} (KL Divergence). The KL divergence between two pdfs $f_1$ and $f_0$ is defined as
$
D (f_1\| f_0) = \int  f_1 (x) \;   \log (f_1 (x)/f_0(x)) \; dx.
$  
\end{definition}
Note that $D (f_1\| f_0) \geq 0$ with equality if and only if $f_1 = f_0$ almost surely. It is usually assumed that
$0 < D(f_1\| f_0) < \infty.$ 

Define the log-likelihood ratio for an observation $X$:
\begin{equation} \label{eq:LLR}
\ell(X) := \log f_1(X)/f_0(X).
\end{equation}
A fundamental property of the log-likelihood ratio, 
which is useful for constructing sequential change detection algorithms, is that before the change $n< \gamma$, the expected value of $\ell(X_n)$ is equal to $-D(f_0||f_1)< 0$; and after the change, $n\geq \gamma$, the expected value of $\ell(X_n)$ is equal to $D(f_1||f_0) >0$. As will be seen later, the KL divergence between the pre- and post-change distributions is an important quantity that characterizes the tradeoff between the average detection delay and the false-alarm rate.

\subsection{Common Sequential Change Detection Procedures}

We now present several commonly used sequential change detection procedures, including the Shewhart chart, CUSUM, and Shiryaev-Roberts procedure, which enjoy certain optimality properties that we will make more precise later in Section~\ref{sec:optimality}. These algorithms can be efficiently implemented in an online setting, which makes them useful in practice. We also briefly discuss some other sequential change detection procedures. 

\subsubsection{Shewhart Chart}

One of the earliest sequential change detection procedures is the Shewhart chart \cite{shew-jamstaa-1925,shew-book-1931}, which is widely used in industrial quality control \cite{montgomery2007introduction}. The Shewhart chart was first introduced for the Gaussian model and based on comparing the instant observation to a threshold. We consider the log-likelihood-based modification and generalization of the standard Shewhart chart, where we compute the log-likelihood ratio based on the current observation (or the current batch of observations) and compare it with a threshold (called the control limit) to make a decision about the change. The property of the log-likelihood ratio discussed in Section \ref{sec:preliminaries} is utilized, which motivates the Shewhart chart:
\[
\taush = \inf \left \{n\geq 1:  \ell(X_n) >   b \right\},
\]
where $b$ is a pre-specified threshold.
The Shewhart chart is widely used in practice due to its simplicity. 
In \cite{pollak2013shewhart}, Pollak and Krieger showed that the Shewhart chart enjoys the optimality property that it maximizes the probability of detecting the change at the time it occurs, subject to false-alarm constraints. 
However, the Shewhart chart may suffer from ``information loss'' due to the fact that it ignores past observations in making a decision about the change, which leads to performance loss when we consider criteria, e.g., the tradeoff between detection delay and false-alarm rate (see Section~\ref{sec:optimality}).

\subsubsection{Cumulative Sum (CUSUM) Procedure}\label{sec:CUSUM}

The CUSUM procedure, first introduced by Page \cite{page-biometrica-1954}, addresses the problem of ``information loss'' in the Shewhart chart. The  CUSUM  procedure uses past observations and thus can achieve a significant performance gain, especially when the change is small. Although the CUSUM procedure was developed heuristically, it was later shown in
\cite{Lorden1971,mous-astat-1986,ritov-astat-1990,lai-ieeetit-1998} that it has very strong optimality properties, which we will discuss further in Section \ref{opt:minimax}. 

The CUSUM procedure utilizes the properties of the cumulative log-likelihood ratio sequence: 
\[
S_n = \sum_{k=1}^n \ell (X_k).
\]
Before the change occurs, the statistic has a negative drift because the expected value of $\ell(X_k)$ before the change is negative. After the change, it has a positive drift because the expected value of $\ell(X_k)$ after the change is positive. Thus, $S_n$ roughly attains its minimum at the change-point $\gamma$. The CUSUM procedure is then constructed to detect this change in the drift of $S_n$.  Specifically, the exceedance of $S_n$ with respect to its past minimum is taken and compared with a threshold $b>0$:
\begin{equation}\label{eq:PageCUSUM_driftbased}
 \tauc = \inf\left\{n\geq 1: W_n = \left(S_n - \min_{0\leq k \leq n } S_k\right)   \geq b \right\}.
\end{equation}
The CUSUM statistic can be rewritten as:
\begin{flalign}\label{eq:wn}
W_n
= \max_{0 \leq k \leq n }\sum_{i=k+1}^n \ell (X_i)
= \max_{1\leq k \leq n+1 }\sum_{i=k}^n \ell (X_i).
\end{flalign}
Note that the maximization over all possible $\gamma = k$ corresponds to plugging in a maximum likelihood estimate of the unknown change-point location in the log-likelihood ratio of the observations to form the CUSUM statistic. 
It can be shown that $W_n$ can be computed recursively:
\[
W_n = (W_{n-1} + \ell (X_n))^+, \quad W_0 = 0,
\]
where $(x)^+ = \max\{x, 0\}$. This recursion enables the efficient online implementation of the CUSUM procedure in practice. 

\subsubsection{Shiryaev-Roberts Procedure}

The maximum likelihood interpretation of the CUSUM procedure is  closely related to another popular algorithm in the literature, called the Shiryaev-Roberts (SR) procedure. In the SR procedure, the maximum in (\ref{eq:wn}) is replaced by a sum and the log-likelihood ratio is replaced by likelihood ratio. The detection statistic for the SR procedure is then defined as:
\begin{equation}\label{SR}
T_n := \sum_{1\leq k \leq n }\prod_{i=k}^n e^{\ell (X_i)},
\end{equation}
and the corresponding stopping time is defined as
\[
\tausr = \inf \{n\geq 1: T_n \geq b\}.
\]
The SR statistic can also be computed recursively:
\[
T_n = (1+T_{n-1}) e^{\ell(X_n)}, \quad T_0 = 0.
\]

\subsection{Optimality}\label{sec:optimality}

We now briefly summarize optimality results in the existing literature for the above procedures. We begin by considering the non-Bayesian setting, where we do not assume a prior on the change-point $\gamma$, and then consider the Bayesian setting, where the change-point is assumed to follow a certain distribution. 

A fundamental problem in sequential change detection is to optimize the tradeoff between the false-alarm rate and the average detection delay, as illustrated in Section \ref{sec:pro_def} using the example in Fig.\,\ref{fig:change-point-example}. Controlling the false-alarm rate is commonly achieved by setting an appropriate threshold on a test statistic such as the one in \eqref{eq:PageCUSUM_driftbased}. But the threshold also affects the average detection delay. A larger threshold incurs fewer false alarms but leads to a larger detection delay, and vice versa. 

\subsubsection{Minimax Optimality}\label{opt:minimax}

In non-Bayesian settings, the change-point is assumed to be a {\it deterministic} and unknown variable. In this case, the average run length ($\ARL$) to false alarm is generally used as a performance measure for false alarms:
\begin{flalign}\label{eq:ARLDef}
\ARL(\tau)=\Expect_\infty[\tau],
\end{flalign}
where $\mathbb P_\infty$ is the probability measure on the sequence of observations when the change never occurs, and $\Expect_\infty$ is the corresponding expectation. Its reciprocal, the false-alarm rate ($\FAR$), is also commonly used:
\begin{equation}\label{eq:FARDef}
\FAR(\tau) = \frac{1}{\ARL(\tau)}=\frac{1}{\Expect_\infty[\tau]}.
\end{equation}
$\FAR$ can also be interpreted as the rate at which false alarms occur in the pre-change regime if we repeat the change detection procedure after each false alarm. Denote the set of stopping times that satisfy a constraint $\alpha$ on the $\FAR$ by:
\begin{equation}\label{eq:DalphaDef}
\Dd_{\alpha}=\{\tau: \FAR(\tau) \leq \alpha\}.
\end{equation}

Finding a uniformly powerful test that minimizes the delay over all possible values of the change-point $\gamma$, subject to a $\FAR$ constraint, is generally intractable. Therefore, it is more tractable to pose the problem in the so-called minimax setting. There are two essential measures of the detection delay in the minimax setting, due to Lorden \cite{Lorden1971} and Pollak \cite{poll-astat-1985}, respectively. 

Lorden considers the supremum of the average detection delay conditioned on the worst possible realizations. In particular, Lorden defines\footnote{Lorden defined $\WADD$ with $(\tau -n +1)^{+}$ inside the expectation, i.e., he assumed a penalty of 1 if the algorithm stops at the change-point. We drop this additional penalty in our definition to be consistent with the other delay definitions in this paper.}:
\begin{equation}\label{eq:WADDdef}
\WADD(\tau) = \underset{n \geq 1}{\operatorname{\sup}}\ \esssup \ \Expect_n\left[(\tau-n)^+| X_1, \dots, X_{n-1}\right], 
\end{equation}
where $\mathbb P_n$ denotes the probability measure on the observations when the change occurs at time $n$, and $\Expect_n$ denotes the corresponding expectation. We then have the following Lorden's formulation:
\begin{equation}
\mbox{minimize } \WADD(\tau) \mbox{ subject to } \FAR(\tau) \leq \alpha. \label{Lorden}
\end{equation}
For the i.i.d. setting, Lorden showed that Page's CUSUM procedure given in \eqref{eq:PageCUSUM_driftbased} is asymptotically optimal as $\alpha \rightarrow 0$. It was later shown in \cite{mous-astat-1986} and \cite{ritov-astat-1990} that a slight modification of the CUSUM procedure, with $W_n = (W_{n-1})^+ + \ell (X_n)$, is exactly optimal for (\ref{Lorden}) for all $\alpha>0$.

Although the CUSUM procedure is exactly optimal under Lorden's formulation, $\WADD(\tau)$ is a pessimistic measure of detection delay since it considers the worst-case pre-change samples. 
An alternative  measure of detection delay was suggested by Pollak \cite{poll-astat-1985}:
\begin{equation}\label{eq:PollakCADDDef}
\CADD(\tau) = \underset{n \geq 1}{\operatorname{\sup}}\ \Expect_n[\tau-n| \tau\geq n], 
\end{equation}
for all stopping times $\tau$ for which the expectation is well-defined. It is easy to see that for any stopping time $\tau$, $\WADD(\tau) \geq \CADD(\tau)$, and therefore, Pollak's formulation is less pessimistic. 

In general, it may be challenging to exactly solve the problem in \eqref{Lorden} and the corresponding problem defined using $\CADD$ in \eqref{eq:PollakCADDDef}. For this reason, asymptotically optimal solutions for the above problems are often investigated in the literature. Specifically, a stopping time $\tau$ is said to be \emph{first-order} asymptotically optimal if it satisfies:
\[
\frac{\CADD(\tau)}{\inf_{\tau \in \Dd_{\alpha}} \CADD (\tau)} \rightarrow 1, \ \ \text{ as } \alpha\rightarrow 0;
\]
it is \emph{second-order} asymptotically optimal if $\CADD(\tau)$  is within a constant of the best possible delay over the class $\Dd_{\alpha}$:
\[{\CADD(\tau)}-{\inf_{\tau \in \Dd_{\alpha}} \CADD (\tau)} =O(1);\] 
 and it is \emph{third-order} asymptotically optimal if such a constant goes to $0$ as $\alpha\to0$:
\[{\CADD(\tau)}-{\inf_{\tau \in \Dd_{\alpha}} \CADD (\tau)} =o(1).
\]
These notions can also be defined similarly for the problem in \eqref{Lorden} defined using $\WADD$.

Pollak's formulation has been studied for the i.i.d. data in \cite{poll-astat-1985} and
\cite{tart-poll-polu-arxiv-2011}. The first-order asymptotic optimality for Lorden's formulation can also be extended to Pollak's formulation. To show this, Lorden in \cite{Lorden1971} established a universal lower bound for $\WADD$ and Lai in \cite{lai-ieeetit-1998} proved the lower bound to $\CADD$:
\begin{theorem}[Lower Bound for $\CADD$ \cite{lai-ieeetit-1998}] 
As $\alpha\rightarrow 0$,
\[
\inf_{\tau \in \Dd_{\alpha}} \CADD (\tau) \geq \frac{|\log \alpha|}{D(f_1||f_0)} (1+o(1)). 
\]
\end{theorem}
It can be shown that the CUSUM procedure with a threshold $b=|\log \alpha|$ is first-order asymptotically optimum for both Lorden's and Pollak's formulations. In particular, as $\alpha\to0$,
\[
\CADD (\tauc) = \WADD (\tauc) \sim  \frac{|\log \alpha|}{D(f_1 || f_0)},
\]
where $\sim$ means the ratio of the quantities on its two sides approaches 1 as $\alpha\rightarrow 0$.

The SR procedure is also asymptotically optimal and it was shown in \cite{tart-poll-polu-arxiv-2011} that by setting the threshold $b=1/\alpha$,
\[\CADD(\tausr) = \frac{|\log \alpha|}{D(f_1||f_0)} +  \xi + o(1),\]
where $\xi$ is a constant that can be characterized using the nonlinear renewal theory \cite{wood-nonlin-ren-th-book-1982} (details omitted here). 

Finally, results in \cite{pollak2013shewhart,moustakides2014multiple,tartakovsky2019sequential} show that the Shewhart chart is optimal for the criterion of maximizing the probability of detecting the change upon its occurrence subject to the $\FAR$ constraints. A more precise statement of this optimality property is as follows. Let the post-change density be denoted by $f_\theta(x)$, where $\theta\in\Theta$ is the post-change parameter. The Shewhart chart as defined earlier becomes the following stopping time:
\[
\taush = \inf \left\{n\geq 1: \frac{f_\theta (X_n)}{f_0(X_n)} > b\right\},
\]
where $b$ is a pre-specified threshold. It is shown that when the threshold $b$ is selected such that $\FAR(\taush)=\alpha$, then $\taush$ is the optimal solution to the following optimization problem:
\begin{equation}
\mbox{maximize } \inf_{1\leq n <\infty} \mathbb P_n^\theta (\tau=n|\tau \geq n)\mbox{ subject to } \FAR(\tau) \leq \alpha, \label{Shewhart}
\end{equation}
where $\mathbb P_n^\theta$ denotes the probability when the change happens at $n$ with $\theta$ being the post-change parameter. Moreover, it was shown in \cite{tartakovsky2019sequential} that if the likelihood ratio $f_\theta (X)/f_0(X)$ is a monotone non-decreasing function of a statistic $S(X)$, then the Shewhart chart is equivalent to $\taush=\inf \{n\geq 1: S(X_n)> b\}$ and when $b$ is selected such that $\FAR(\taush)=\alpha$, the Shewhart chart is uniform optimal in $\theta\in\Theta$ in the sense of solving \eqref{Shewhart} for all $\theta\in\Theta$.

\begin{figure}[!ht]
\center
\includegraphics[width=0.8\linewidth]{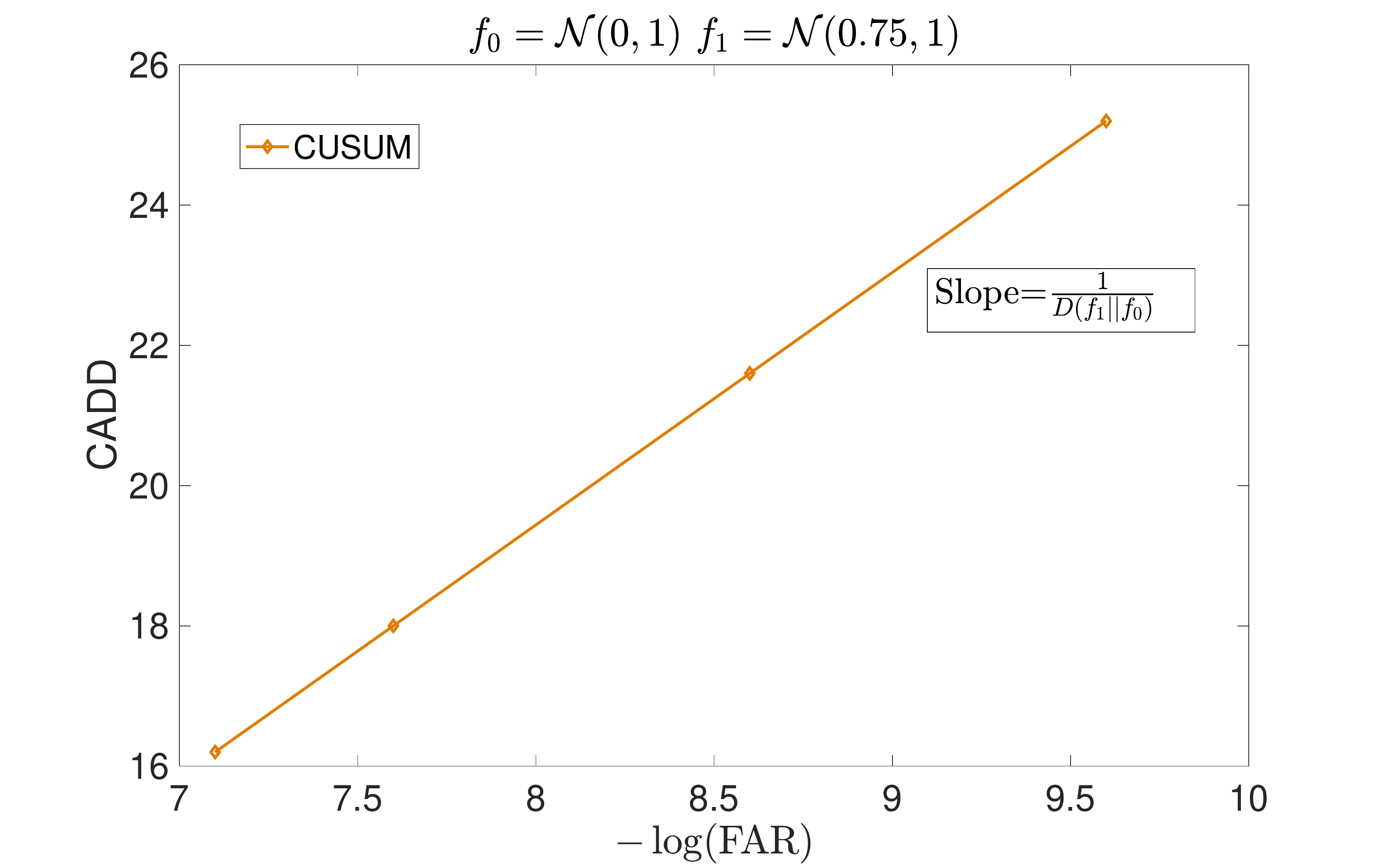}
\caption{Tradeoff curve between $\CADD$ and $-\log(\FAR)$ for the CUSUM algorithm. The pre-change distribution is $f_0=\mathcal{N}(0,1)$, and the post-change distribution is $f_1=\mathcal{N}(0.75,1)$. The slope of the curve is approximately $1/D(f_1 || f_0)$.}
\label{fig:CUSUM_Tradeoff_slope}
\end{figure}

In summary, both the CUSUM and SR procedures are asymptotically optimal with respect to Lorden's formulation and Pollak's formulation. The FAR decays to zero exponentially with exponent $D(f_1||f_0)$. We demonstrate the theory using an example in Fig.\,\ref{fig:CUSUM_Tradeoff_slope} by plotting the tradeoff curve between the $\CADD$ and $-\log(\FAR)$ for the CUSUM procedure. Note that the curve has a slope approximately of $1/D(f_1 || f_0)$, which is consistent with the theory. 
 
Some more optimality results are summarized as follows.  Under Pollak's criterion, it was shown in \cite{tart-poll-polu-arxiv-2011} that the SR algorithm is second-order asymptotically optimal, and that the SRP algorithm (Pollak's version of the SR algorithm that starts from a quasi-stationary distribution of the SR statistic) is third-order asymptotically optimal (as was also first established in \cite{poll-astat-1985}). More importantly, in \cite{tart-poll-polu-arxiv-2011}, it was proved that the SR-$r$ procedure that starts from a specially selected fixed point $r$ is third-order optimal. In \cite{polunchenko2010optimality}, it was shown that SR-$r$ is strictly optimal for $\CADD$ in some special cases. We also note that the (generalized) Shewhart chart is optimal for the criterion of maximizing the probability of detecting a change subject to false alarm constraints.

\subsubsection{Bayesian Optimality}\label{sec:bayes_optimaly}

In the Bayesian setting, it is assumed that the change-point
is a random variable $\Gamma$ taking values on the non-negative integers,
with probability mass function $\pi_n = \Prob\{\Gamma =n\}$. For a stopping time $\tau$, define the average detection delay ($\ADD$) and the probability of false alarm ($\PFA$) as follows:
\begin{eqnarray}
\ADD(\tau)&=&\Expect \left[(\tau - \Gamma)^+\right] = \sum_{n=0}^\infty \pi_n \Expect_n   \left[(\tau - \Gamma)^+\right], \label{eq:ADD_def}\\
\PFA(\tau) &=& \Prob (\tau<\Gamma) = \sum_{n=0}^\infty \pi_n \Prob_n (\tau<\Gamma).
\end{eqnarray}
In Bayesian sequential change detection, the goal is to minimize $\ADD$  subject to a constraint on $\PFA$. 
Shiryaev \cite{shir-siamtpa-1963} formulated the Bayesian sequential change detection problem   as follows:
\begin{equation}\label{eq:BayesConstProb}
\mbox{minimize } \ADD(\tau) \mbox{ subject to } \PFA(\tau) \leq \alpha. \quad (\mbox{Shiryaev})
\end{equation}
The prior on the change-point $\Gamma$ is usually assumed to be a geometric distribution with parameter $0 < \rho < 1$, 
\begin{equation} \label{eq:geom_prior}
\pi_n  = \Prob \{ \Gamma= n \} =  \rho (1-\rho)^{n-1}\mathbb I_{\{n \geq 1\}}, \quad  \pi_0 =0,
\end{equation}
where $\mathbb I$ is the indicator function. 
The justification for this assumption is that the geometric distribution is memoryless. Moreover, it leads to a tractable formulation and convenient optimal solutions to the Bayesian problem in \eqref{eq:BayesConstProb} as we will discuss in the following.

The detection statistic of the \emph{Shiryaev algorithm}  is the posterior probability that the change has taken place given the observations so far. Denote by $X_1^n = (X_1, \ldots, X_n)$ the observations up to time $n$, and by
\begin{equation} \label{eq:pk}
p_n = \Prob ( \Gamma \leq n \; | \; X_1^n )
\end{equation}
the \textit{a posteriori} probability at time $n$ that the change has taken place given the observations up to time $n$. It then follows from the Bayes' rule that $p_n$ can be updated recursively:
\begin{equation} \label{eq:pkrec}
p_{n+1} = \frac{\tilde{p}_n e^{\ell(X_{n+1})}}{ \tilde{p}_n  e^{\ell(X_{n+1})} + (1-\tilde{p}_n )},
\end{equation}
where $\tilde{p}_n = p_n + (1-p_n) \rho$,
and $p_0=0$.
Then the \emph{Shiryaev algorithm} is defined by comparing $p_{n}$ with a given threshold $b_\alpha$:
\begin{equation}
\label{eq:OptimalAlgo}
\taus = \inf\left\{ n\geq 1: p_n \geq b_{\alpha} \right\},
\end{equation}
where $b_\alpha \in (0,1)$ is chosen such that the false alarm constraint, $\PFA (\taus)  \leq \alpha$, is satisfied.
\begin{theorem}[Optimal Bayesian Procedure \cite{shir-siamtpa-1963,shir-opt-stop-book-1978}]\label{thm:ShiryaevOpt} When the threshold $b_\alpha$ is selected such that $\PFA(\tau_{\scriptscriptstyle \text{\emph{S}}}) = \alpha$, the Shiryaev algorithm in \eqref{eq:OptimalAlgo} is Bayesian optimal for \eqref{eq:BayesConstProb}.
\end{theorem}

An equivalent form of the Shiryaev statistic can be developed using the idea of the likelihood ratio test. This builds a connection to the earlier SR statistic defined in (\ref{SR}), and it reveals useful insights about the nature of the procedure. Consider two hypotheses: ``$H_1: \Gamma \leq n$'' and ``$H_0: \Gamma > n$''. Denote by $R_{n,\rho} = p_n/[\rho(1-p_n)]$ the scaled likelihood ratio between the two hypotheses averaged over the change-point. 
It then follows that $R_{n,\rho}$ can be updated recursively as:
\begin{equation}\label{R_recursion}
R_{n+1,\rho} = \frac{1+R_{n,\rho}}{1-\rho} e^{\ell (X_{n+1})}, \quad R_{0,\rho} = 0.
\end{equation}
The Shiryaev stopping time $\taus$ in \eqref{eq:OptimalAlgo} can then be rewritten as a comparison of  $R_{n,\rho}$ with a threshold.
We remark here that if we set $\rho = 0$, then the Shiryaev statistic reduces to the SR statistic in (\ref{SR}).

A generalized Shewhart chart is also Bayesian optimal, as shown in \cite{pollak2013shewhart}, in the sense that it minimizes the expected loss where the loss function is $ \mathbb I_{\{\tau \neq \Gamma\}}$, assuming that the change-point $\Gamma$ follows a geometric prior and the parameter $\theta$ of the post-change distribution follows a known prior $G$. This result was generalized in  \cite[Theorem 5.1]{tartakovsky2019sequential}. Moreover, both the CUSUM and SR procedures are first-order asymptotically optimal for the Bayesian setting when the prior has a heavy tail, or when the change-point is geometrically distributed with a small enough parameter.

\subsubsection{Evaluating the Performance Metrics} \label{sec:eval}

In the definition of the $\WADD$ metric \eqref{eq:WADDdef} and the $\CADD$ metric \eqref{eq:PollakCADDDef}, it appears that we need to consider the supremum over all possible past observations and over all possible change-points. However, we can actually show that for the CUSUM and SR procedures, and some other algorithms, that the supremum over all possible change-points in $\WADD$ and $\CADD$ is achieved at time $n=1$:
\[
\begin{aligned}
\CADD (\tauc) &= \WADD (\tauc) = \Expect_1 \left[ \tauc - 1\right],\\
\CADD (\tausr) &= \WADD (\tausr) = \Expect_1 \left[ \tausr - 1\right].
\end{aligned}
\]
Therefore, the  $\CADD$ and the $\WADD$ can be conveniently evaluated by setting $\gamma=1$, without ``taking the supremum''.

\subsection{Other Sequential Change Detection Procedures} \label{sec:other_seq_ch}

\subsubsection{Mixture and Generalized Likelihood Ratio (GLR) Statistics}

The CUSUM and SR procedures require full knowledge of pre- and post-change distributions to obtain the log-likelihood ratio $\ell(X)$ used in computing the test statistics. In practice, the post-change distribution $f_1$ may be unknown. In the parametric setting, the post-change distribution can be parametrized using $f_\theta$, where $\theta \in \Theta$ is the unknown parameter. Two commonly used methods for the situation here, which corresponds to the problem of composite hypothesis testing, are the generalized likelihood ratio (GLR) approach and the mixture approach. 
In the GLR approach, a supremum over $\theta \in \Theta$ is taken in constructing the test statistic. In particular, the test statistic for the GLR-CUSUM algorithm is given by:
\begin{flalign}\label{eq:wnGLR}
   W_n^{\text{G}}
= \max_{1\leq k \leq n+1 } \sup_{\theta \in \Theta} \sum_{i=k}^n \ell_\theta (X_i),
\end{flalign}
where $\ell_\theta (X) = \log (f_\theta(X)/f_0 (X))$. Performance analyses of the GLR-CUSUM algorithm for one-parameter exponential families can be found in \cite{Lorden1971,lorden1973open}. A major drawback of the GLR approach is that the corresponding GLR statistic (e.g., the one given in \eqref{eq:wnGLR}) cannot be computed recursively in time, except in some special cases (e.g., when the parameter set $\Theta$ has finite cardinality). To reduce the computational cost, a window-limited GLR  approach was developed in \cite{willsky1976generalized} and generalized in \cite{lai-ieeetit-1998,lai1999efficient}. Window-limited versions of the GLR algorithm can be shown to be asymptotically optimal in certain cases if the window size is carefully chosen as a function of $\FAR$. 

The mixture method replaces the supremum over $\theta \in \Theta$ by a weighted average. For example, the mixture-CUSUM statistic is computed as:
\begin{flalign}\label{eq:mGLR}
   W_n^{\text{m}}
= \max_{1\leq k \leq n+1 }  \log\int_{\Theta}\prod_{i=k}^n\frac{f_\theta (X_i)}{f_0(X_i)} \omega (\theta) d\theta,
\end{flalign}
where $\omega(\theta)$ is a weight function that integrates (sums) to 1 over $\Theta$.  Note that, like the GLR test statistic, the mixture test statistic cannot be computed recursively in general. It was shown in \cite{tartakovsky2019sequential} that the mixture approach can result in  first-order asymptotically optimal tests for practically any prior for both the i.i.d. and non-i.i.d. cases. In \cite{siegmund2008minimax}, the optimal prior was established such that the resultant mixture SR procedure is asymptotically optimal in a certain stronger sense.

\subsubsection{EWMA}
Note that the CUSUM and SR procedures can achieve a significant gain in performance when compared to the Shewhart chart by making use of past observations, i.e., CUSUM and SR have memory. 
The exponentially weighted moving average (EWMA) chart is another type of sequential change detection procedure that employs past observations. The EWMA detection statistic was originally defined as $Z_n = \lambda X_n + (1-\lambda) Z_{n-1}$, where $\lambda \in (0,1]$ is a pre-specified constant, with the aim to detect mean shift. The EWMA can be generalized to $Z_n = \lambda \ell(X_n) + (1-\lambda) Z_{n-1}$ to detect shift in distribution, more generally. Thus, $Z_n$ is a weighted moving average of all past information with weights decreasing exponentially in time. 
The EWMA chart is simple to implement and does not require any prior knowledge of the pre- and post-change distributions. A performance comparison of the EWMA chart and the CUSUM and SR procedures is given in \cite{polunchenko2014optimal}.

\section{Generalizations and extensions}\label{sec:new-advances}

\subsection{General Asymptotic Theory for Non-i.i.d. Data}
\label{sec:LAI}

There has been a considerable amount of effort to generalize the optimality results for sequential change detection to the non-i.i.d.\ setting. Lai \cite{lai-ieeetit-1998} initiated the development of  a general minimax asymptotic theory for both Lorden's and Pollak's formulations, while  Tartakovsky and Veeravalli \cite{tartakovsky2005general} initiated the development of a general Bayesian asymptotic theory.

\subsubsection{General Minimax Asymptotic Theory}
Under the minimax setting, Lai in \cite{lai-ieeetit-1998} obtained a general lower bound for non-i.i.d. data on the $\CADD$ (and hence on the $\WADD$) for any stopping time that satisfies the constraint that $\FAR$ is no larger than $\alpha$. It was then shown that an extension of the CUSUM procedure \eqref{eq:PageCUSUM_driftbased} to the non-i.i.d.\ setting achieves this lower bound asymptotically as $\alpha\to 0$. There are also works investigating non-i.i.d. data under some specific settings, e.g., multi-sensor slope change detection  \cite{cao2018multi}, linear regression models \cite{Geng2019Linear,wang2019statistically}, generalized autoregressive conditional heteroskedasticity (GARCH) models \cite{berkes2004sequential}, non-stationary time series \cite{choi2008sequential}, general stochastic models  \cite{Tartakovsky2017noniid,tartakovsky2019asymptotically}, and hidden Markov models \cite{Tartakovsky2019Markov}. We refer to \cite{tartakovsky2019sequential} for more recent developments on this topic.

We now present a generalized CUSUM procedure for non-i.i.d.\ data. In this setting, conditional distributions are used to compute the likelihood ratios. In the pre- and post-change regimes, the conditional distribution of $X_n$ given $X_1^{n-1}$ is given by $f_{0,n}(X_n|X_1^{n-1})$ and  $f_{1,n}(X_n|X_1^{n-1})$, respectively. Define the conditional log-likelihood ratio and the CUSUM statistic, respectively, as:
\[
Y_i = \log \frac{f_{1,i}(X_i|X_1^{i-1})}{f_{0,i}(X_i|X_1^{i-1})},\quad \text{and  }
C_n = \max_{1\leq k\leq n+1} \sum_{i=k}^n Y_i.
\]
Then the stopping time for the generalized CUSUM is defined as:
\begin{equation}
\label{eq:CUSUMgen}\taug = \inf\left\{ n\geq 1: C_n \geq b \right\}.
\end{equation}
Note the generalized CUSUM (for non-i.i.d. data) takes a similar form as the original CUSUM (for i.i.d. data) except that we replace the log-likelihood ratio with the conditional log-likelihood ratio. 

The minimax optimality of the generalized CUSUM for the non-i.i.d. data was established in \cite{lai-ieeetit-1998}. Under some regularity conditions, by setting the threshold $b=|\log \alpha|$,  we have $\taug\in\Dd_{\alpha}$. If there exists $I$ such that
\begin{equation}
\label{eq:noniidCond_exist_q}
\underset{m\leq t }{\operatorname{\max}} \frac{1}{t}\sum_{i=n}^{n+m} Y_i \to I \ \ \ \text{ a.s. }\Prob_n, \ \ \text{ as } t \to \infty  \ \ \forall n,
\end{equation}
and the convergence is complete in the sense that $\sum_{n=1}^\infty\mathbb P_1(|(1/n) \sum_{i=1}^n Y_i - I| \geq \epsilon) < \infty$ for all $\epsilon>0$, then as $\alpha\to 0$:
\begin{equation} \label{eq:Lai_lb}
\begin{aligned}
\CADD(\taug) &\sim \WADD(\taug) \sim  \underset{\tau\in \Dd_{\alpha}}{\operatorname{\inf}}\ \WADD(\tau) \\
    & \sim \underset{\tau\in \Dd_{\alpha}}{\operatorname{\inf}}\ \CADD(\tau) \sim \frac{|\log \alpha|}{I},
\end{aligned}
\end{equation}
where the positive constant $I>0$ plays a similar role as the KL divergence in the i.i.d. setting.

\subsubsection{General Bayesian Asymptotic Theory}

Under the Bayesian setting, when the samples conditioned on the change-point are non-i.i.d., it is generally difficult to find an exact solution to the Shiryaev problem in \eqref{eq:BayesConstProb}. Tartakovsky and Veeravalli \cite{tartakovsky2005general} showed that the Shiryaev algorithm is asymptotically optimal as $\alpha\to0$, under some regularity conditions on the pre- and post-change distributions.

Similar to the i.i.d.\ case, we can define the posterior probability $p_n$ of change having occurred before time $n$ given all previous samples, in the same expression as \eqref{eq:pk}. 
The Shiryaev algorithm for the non-i.i.d. setting is then defined in the same way as in \eqref{eq:OptimalAlgo}. 
Note that the recursion in \eqref{eq:pkrec} may not hold for a general distribution for $\Gamma$. However, if the change-point $\Gamma$ is geometrically distributed, a recursive expression for $p_n$ can still be derived. 
Define 
\[
d = - \lim_{n\to \infty} \frac{\log \Prob(\Gamma > n)}{n},
\]
which captures the decay rate of the tail probability of change-point $\Gamma$'s prior distribution as the sample size $n$ increases. 
When $\Gamma$ is ``heavy-tailed'', $d=0$, and when $\Gamma$ has an ``exponential tail'', $d>0$. For example, when the prior distribution is geometric with parameter $\rho$ as defined in \eqref{eq:geom_prior}, $d=|\log(1-\rho)|$.  If there exists $I$ such that
\begin{equation}
\label{eq:noniidCond_exist_q2}
\frac{1}{t}\sum_{i=n}^{n+t} Y_i \to I \ \ \ \text{ a.s. }\Prob_n \ \ \text{ as } t \to \infty,  \ \ \forall n,
\end{equation}
and some additional conditions on the rates of convergence are satisfied (see \cite{tartakovsky2005general} for the details), then the
Shiryaev algorithm in \eqref{eq:OptimalAlgo} with a threshold $b_\alpha=1-\alpha$ is asymptotically optimal for Bayesian optimization problem in \eqref{eq:BayesConstProb} as $\alpha \to 0$ \cite{tartakovsky2005general}:
\begin{equation}
\label{eq:noniid_UB2}
\ADD(\taus)
\sim \inf_{\tau: \PFA(\tau)\leq \alpha}\ADD(\tau)\sim\frac{|\log \alpha|}{I + d}.
\end{equation}
Note that in \cite{tartakovsky2005general}, a general result for the $m$-th moment of the delay was developed. Here, for simplicity, we only presented the result for $m=1$.

\subsection{Change-of-measure to Obtain Accurate $\ARL$ Approximations} 
\label{sec:change-of-measure}

For CUSUM and SR procedures with i.i.d. samples, it may be relatively easy to evaluate their performance (such as the $\ARL$) both theoretically and numerically, as discussed in Section \ref{sec:eval}. However, in many settings such as those involving non-i.i.d. observations, GLR statistics \cite{siegmund2011detecting}, and non-parametric statistics \cite{li2019scan}, it may be challenging to develop exact analytical expressions for the $\ARL$ (or its inverse the $\FAR$). In these situations, one has to use onerous numerical simulation to obtain a threshold for a target $\ARL$. To tackle this problem, techniques based on extremes in random fields have been developed \cite{yakir2013extremes}, from which one can obtain accurate approximations to the $\ARL$ for many problems. 

\subsubsection{Using Change-of-measure to Analyze the $\ARL$}

The main idea here is to relate finding $ \ARL $ to finding the tail probability of the maximum of a random field. To obtain a more accurate approximation of the $\ARL$, an alternative probability measure is considered, under which false alarms are more likely to occur. This is analogous to ``importance sampling'', but it is more involved since the alternative probability measure is usually a mixture of distributions. 

The analysis usually involves two steps. First, we aim to find the probability $\mathbb P_\infty\{\tau\leq m\}$, for a large constant $m>0$ and stopping time $\tau$ (the first time the detection statistic exceeds the threshold $b$). Finding this probability is challenging because $\{\tau\leq m\}$ is a rare event under the pre-change regime, especially when the threshold $b$ is large (this is the asymptotic scenario that we are interested in). Therefore, the change-of-measure technique plays an important role by considering an alternative measure under which $\{\tau\leq m\}$ happens with a much higher probability. More specifically, we choose the alternative measure such that the expectation of the detection statistic equals the threshold $b$. Then, using the local central limit theorem and the local behavior of the correlated random field, we can obtain an analytical expression for the probability of $\{\tau\leq m\}$ under the alternative measure. 
The probability under the alternative measure is then converted back to the probability under the original measure through Mill's ratio. The rigorous mathematical derivations can be found in \cite{yakir2013extremes}.

Second, we will relate the above probability to the $\ARL$, leveraging the fact that the stopping time $\tau$ as threshold $b\rightarrow \infty$ is asymptotically exponentially distributed \cite{aldous2013probability,siegmund2008detecting}. 
Although this fact only holds strictly for stopping times for algorithms such as the CUSUM and SR when observations are i.i.d.  \cite{pollak2009asymptotic}, this method has been widely used and is verified to be highly accurate in practice (see examples in \cite{siegmund2011detecting,xie2013sequential,cao2018multi,li2019scan}). Thus, for a large $m$, 
$
\Prob_\infty \{\tau \leq m\} \sim 1-e^{-\lambda_b  m},
$
where $\lambda_b$ is the parameter of the exponential distribution. By definition, the mean of the exponential distribution is $1/\lambda_b$, which corresponds to the $\ARL$. 

\subsubsection{Example: Analyzing MMD-based Sequential Change Detection Procedure}\label{change-of-measure}

Below, we illustrate the change-of-measure technique by analyzing the non-parametric kernel-based maximum mean discrepancy (MMD) statistics (details can be found in \cite{li2019scan}). 
The kernel MMD divergence, which measures the distance between two arbitrary distributions, is widely adopted in signal processing and machine learning. Given two sets of samples $X:=\{x_1, \ldots, x_n\}$ generated i.i.d. from a distribution $f_0$ and $Y:=\{y_1, \ldots, y_n\}$ generated i.i.d. from a distribution $f_1$, an unbiased estimator of the MMD between $f_0$ and $f_1$ is the following:
\[
\begin{aligned}
\mbox{MMD}(X, Y) = \frac{1}{n(n-1)} \sum_{i\neq j} \{ & k(x_i, x_j) + k(y_i, y_j) \\
&- k(x_i, y_j) - k(x_j, y_i)\},
\end{aligned}
\]
where $k(\cdot, \cdot)$ is a kernel in the reproducing kernel Hilbert space (RKHS), e.g., Gaussian kernel. Intuitively, the MMD statistic is small when $f_0$ is similar to $f_1$, and is large otherwise.

The sequential change detection procedure based on the MMD statistic is then defined as follows \cite{li2019scan}. At each time $t$, we treat the most recent $B$ samples, denoted by $X_{t-B+1}^t:=\{X_{t-B+1}, \ldots, X_t\}$, as the test block ($B > 0$ is a pre-specified parameter). Then we sample $N$ blocks of size $B$ from the ``reference'' data generated from the pre-change distribution, denoted by $\{\tilde X_1,\ldots,\tilde X_N\}$. We compute an average MMD statistic of all the reference blocks with respect to the test block:
\[
U_t = \frac 1 N \sum_{i=1}^N \mbox{MMD}(\tilde X_i, X_{t-B+1}^t). 
\]
Define $Z_t':=U_t/\sqrt{\mbox{var}[U_t]}$ as the standardized detection statistic, where the variance $\mbox{var}[U_t]$ can be found in closed-form and can be estimated conveniently from data \cite{li2019scan}. 
The MMD-based procedure stops when the standardized MMD statistic exceeds a threshold $b$:
\[
\taum = \inf \{t: Z_t'> b\}.
\]
This corresponds to a generalized type of Shewhart chart.
\begin{theorem}[$\ARL$ of MMD-based Procedure \cite{li2019scan}]\label{ARL_MMD}
Let $B > 0$. When $b\rightarrow \infty$, the $\ARL$ of the stopping time $\tau_{\scriptscriptstyle \text{\emph{M}}}$, $\mathbb E_\infty[\tau_{\scriptscriptstyle \text{\emph{M}}}]$, is given by:
\[
\frac{e^{b^2/2}}{b} \left\{\frac{2B-1}{\sqrt{2\pi} B(B-1)} \nu\left(b \sqrt{\frac{2(2B-1)}{B(B-1)}}\right)\right\}^{-1}(1+o(1)),
\]
where $\nu(\cdot)$ is a special function whose definition can be found in \cite{Siegmund1985}. 
\end{theorem}
We present the main step of the proof to Theorem\,\ref{ARL_MMD} to illustrate the change-of-measure technique. First, note that the event $\{\tau \leq m\}$ is the same as the maximum of the detection statistic has exceed the threshold $b$ at some point before $m$, i.e., $\{\sup_{2\leq t\leq m}Z_t' \geq b\}$, and
\begin{align*}
&\mathbb P_\infty\left\{\sup_{2\leq t\leq m}Z_t' \geq b\right\} \\
& = \mathbb E_\infty\left\{ \frac{\sum_{t=2}^m e^{\xi_{t}}}{\sum_{s=2}^me^{\xi_s}}\mathbb I_{\{ \sup_{2\leq t\leq m}Z_t' \geq b \}} \right\} \\
&= e^{-b^2/2}\sum_{t=2}^m \mathbb E_t \big\{R_t e^{- [\xi_t -b^2/2+ \log M_t] }\mathbb I_{\{\xi_t - b^2/2+ \log M_t \geq 0\}}\big\},
\end{align*}
where $\xi_t = bZ_t' - b^2/2$ is the log-likelihood ratio between the changed measure $\mathbb E_t[X] = \mathbb E_\infty[X e^{\xi_t}]$ and the original measure $\mathbb E_\infty$, $M_t = \max_s e^{\xi_s - \xi_t}$, $S_t = \sum_s e^{\xi_s - \xi_t}$, and the so-called {\it Mill's ratio} $R_t = M_t/S_t$. The result in Theorem\,\ref{ARL_MMD} is established by establishing properties of the local field $\{\xi_s - \xi_t\}$ and the global term $\xi_t - b^2/2$ (details omitted here).

\begin{figure}[!ht]
    \centering
    \includegraphics[width = 0.7\linewidth]{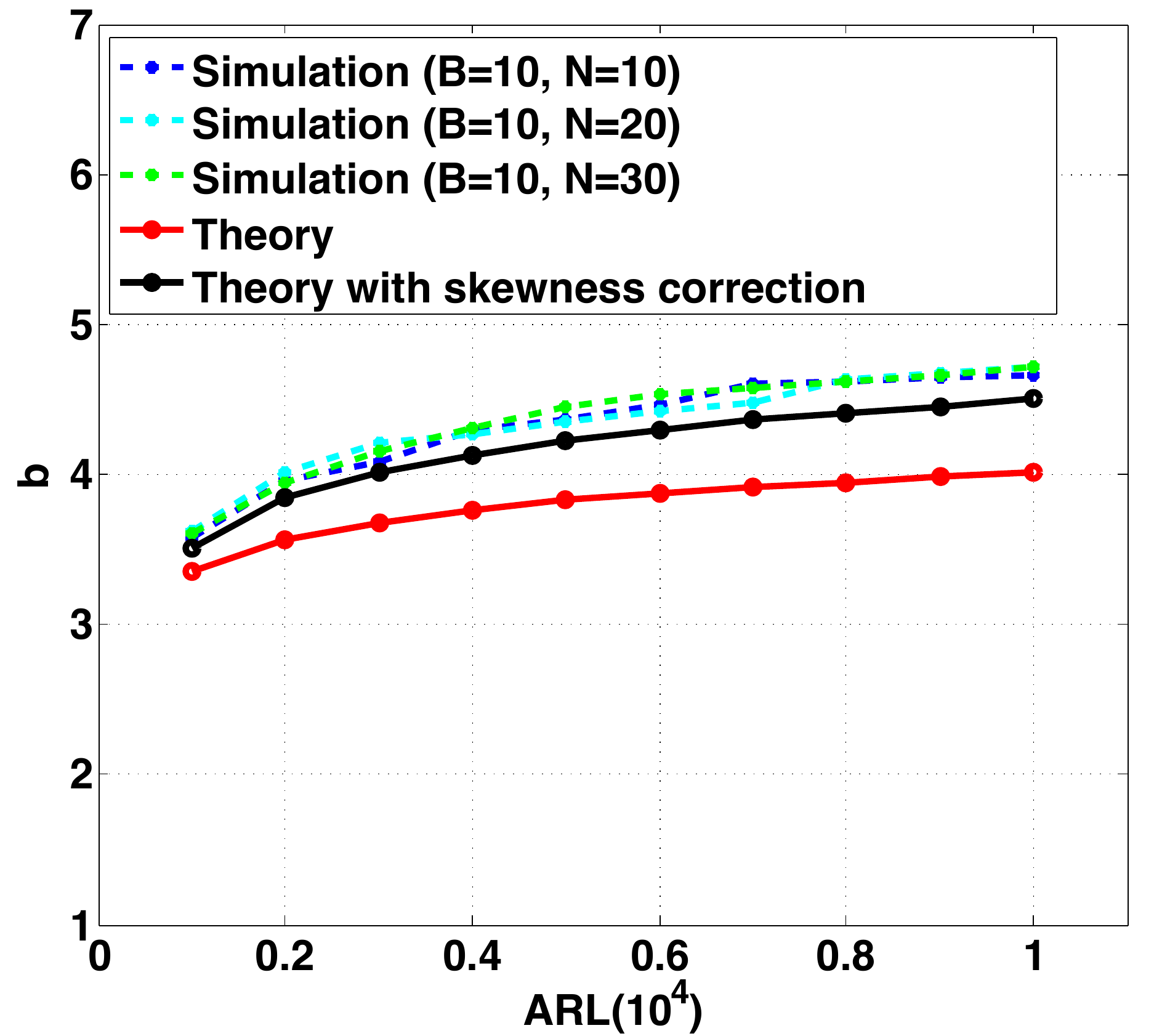}
    \caption{Accuracy of $\ARL$ approximations, obtained by ``change-of-measure'', for the sequential MMD-based procedure: comparison of the thresholds obtained by simulation and from Theorem~\ref{ARL_MMD}.} 
    \label{fig:ARLN}
\end{figure}

The numerical example in Fig.\,\ref{fig:ARLN} demonstrates that the threshold $b$ (to achieve a target $\ARL$) obtained using the theoretical approximation in Theorem\,\ref{ARL_MMD} is consistent with that obtained from simulations, especially after a skewness correction. This example demonstrates that the theoretical approximation of the $\ARL$ obtained using the change-of-measure technique is of high accuracy, and thus can help avoid computationally expensive simulations to calibrate the procedure. 

\subsection{Non-stationary and Multiple Changes}

In various modern applications, for instance, line outage detection in power systems \cite{Rovatsos2:2016} and stochastic power supply control in data centers \cite{rovatsos2018asilomar}, the change is not stationary. There can be a sequence of multiple changes: one followed by another. Below, we review some recent advances in sequential detection of dynamic changes.

\subsubsection{Sequential Change Detection Under Transient Dynamics}\label{sec:transient}

In classical sequential change detection formulations \cite{veeravalli2013quickest,detectAbruptChange93,tartakovsky2014sequential,poor-hadj-QCD-book-2008}, the statistical behavior of the observations is characterized by one pre-change distribution and \textit{one} post-change distribution (known or unknown). In other words, the statistical behavior after the change is stationary. This assumption may be too restrictive for many practical applications with more involved statistical behavior after the change-point.

An example of the problem where the observations are non-stationary after the change, is sequential change detection under transient dynamics, which was studied in \cite{Zou2017QCD,rovatsos2018asilomar,george2017icassp,Rovatsos2:2016}. Specifically, the pre-change distribution does not change to a persistent  post-change distribution instantaneously, but after several transient phases, each phase is associated with a distinct data generating distribution. The goal is to detect the change as quickly as possible, either during the transient phases or during the persistent phase. This problem is fundamentally different from detecting a transient change (see, e.g., \cite{guepie2012,ebrahimzadeh2015sequential,egea2018performance}), where the system goes back to the pre-change mode after a single transient phase, and the goal is to detect the change within the transient phase. The problem is also related to sequential change detection in the presence of a nuisance change, where the presence of the nuisance change can be modeled as a transient phase. However, an alarm should be raised only if the critical change occurs \cite{lau2019quickest}.

Two algorithms were proposed and investigated in \cite{Zou2017QCD,Rovatsos2:2016} for the minimax setting, the dynamic-CUSUM (D-CUSUM), and the weighted dynamic-CUSUM (WD-CUSUM), where the change-point and the transient durations are assumed to be unknown and deterministic. The basic idea is to construct a generalized likelihood based algorithm taking the supremum over the unknown change-point and the durations of transient phases. It was shown in \cite{Zou2017QCD,Rovatsos2:2016} that the D-CUSUM and WD-CUSUM test statistics can be updated recursively, and thus are computationally efficient. In \cite{Zou2017QCD}, it was demonstrated that both algorithms are adaptive to the unknown transient dynamics, although durations of transient phases were unknown and were not employed in algorithm implementation. Moreover, both the D-CUSUM (under certain conditions) and the WD-CUSUM algorithms were shown to be first-order asymptotically optimal in \cite{Zou2017QCD}. The Bayesian setting was investigated in \cite{george2017icassp}, where the change-point and the \textit{durations} of transient phases are assumed to be geometrically distributed. The optimal test was constructed, and a computationally efficient alternative test based on thresholding the posterior probability that the change has occurred was also proposed.


\subsubsection{Sequential Detection of Moving Anomaly}
Existing studies on sequential change detection in networks usually assume that the change is persistent once it affects a node. However, there are scenarios where the change may not necessarily be persistent at a particular node; instead, it is persistent across the network as a whole, e.g., a moving anomaly in a sensor network. In this case, existing approaches using CUSUM statistics from each node, e.g., \cite{mei2010efficient,zou2019quickestnetwork,fellouris2017multichannel,hadjiliadis2009one}, cannot be applied. Recently, the problem of sequential moving anomaly detection in networks was studied in \cite{rovatsos2020sequential,rovatsos2020sequentialmixture}. Specifically, after an anomaly emerges in the network, one node is affected by the anomaly at each time instant and receives data from a post-change distribution. The anomaly dynamically moves across the network with an unknown trajectory, and the node that it affects changes with time. Two approaches have been proposed to model the trajectory of the anomaly: the hidden Markov model \cite{rovatsos2020sequential}, and the worst-case approach \cite{rovatsos2020sequentialmixture}, which we discuss in the following. 

The first approach (hidden Markov model) \cite{rovatsos2020sequential} models the anomaly's trajectory as a Markov chain, and thus the samples are generated according to a hidden Markov model. The advantage of this model is that it takes into consideration the network's topology, i.e., that the anomaly only moves from a node to one of its neighbors. In \cite{rovatsos2020sequential}, a windowed GLR based algorithm was constructed and was shown to be first-order asymptotically optimal. Alternative algorithms were also designed with performance guarantees, including the dynamic SR procedure, recursive change-point estimation, and a mixture CUSUM algorithm.

The second approach (worst-case approach) \cite{rovatsos2020sequentialmixture} assumes that the anomaly's trajectory is unknown but deterministic and considers the worst-case performance over all possible trajectories. A CUSUM-type procedure was constructed. The main idea is to use the mixture likelihood to construct a test statistic, which is further used to build a procedure of the CUSUM-type. This procedure was shown to be exactly optimal in \cite{rovatsos2020sequentialmixture} when the sensors are homogeneous. This idea has been further generalized to solve the sequential moving anomaly detection problem with heterogeneous sensors and has been shown to be first-order asymptotically optimal \cite{rova-veer-mous-isit2020}. 

\subsubsection{Multiple Change Detection}

A related line of research is multiple change detection in the offline setting, which aims to estimate multiple change-points from observations in a retrospective study. Various methods were proposed to estimate the number and locations of change-points, including hierarchical clustering based method \cite{Nonpara-multu-CPD2014}, binary segmentation type methods \cite{vostrikova1981detecting,bai1997estimating,Cho2012Multiscale,Fryzlewicz2014BS,Cho2015BS-high-dim,fryzlewicz2014wild}, (penalized) least-squared methods \cite{YaoAu1989,lavielle2000least,lavielle2005using,lebarbier2005detecting,boysen2009consistencies}, Schwarz criterion \cite{yao1988estimating}, kernel-based algorithms \cite{arlot2019kernel,Harchaoui2007}, and so on. Another line of work aims to reduce the computational complexity of the multiple change detection methods, such as \cite{killick2012optimal,rigaill2010pruned,harchaoui2010multiple}. We refer to \cite{TRUONG2020107299} for a recent review on multiple change detection. Some offline multiple change detection algorithms can motivate the development of their online versions.

\subsubsection{Decentralized and Asynchronous Change Detection in Networks} 
When the information for detection is distributed across a network of sensors, detection problems fall under the umbrella of distributed (or decentralized) detection \cite{tsitsiklis1993decentralized,vars-book-1996,cham-veer-ieeespm-2007,veer-vars-survey-2012}. In the decentralized setting, each sensor sends messages to the fusion center based on the observations it has received so far. The fusion center may provide feedback to sensors and make the final decision. The problem of decentralized sequential change detection in distributed sensor systems was introduced in \cite{veeravalli2001decentralized}, considering the observation model where all sensors are affected by the change at the same time. There have been a number of papers on the topic since then, see e.g., \cite{tartakovsky2004change,tartakovsky2008asymptotically,mei2005information}. A more recent (and practical) perspective is that the change may affect sensors with delay, i.e., different sensors may observe the change at different times, which we will present in the following.

In the case of multiple data streams, the change may happen asynchronously for different sensors. When we desire to detect the first onset of change, it is proposed in \cite{hadjiliadis2009one} to monitor each data stream by {\it local} CUSUM procedures and raise the alarm when any sensor raises an alarm. The sum of local CUSUM statistics has been considered in  \cite{mei2010efficient} and was shown to be asymptotically optimal. The problem where the change propagates from one sensor to the next with known Markov dynamics after the change was studied in \cite{VenuSensor2010}, and an asymptotically optimal test was developed. A recent procedure proposed in \cite{xie2019asynchronous} finds an optimal combination of local data streams accounting for their delays in being affected by the change, which can boost the signal-to-noise ratio and reduce the detection delay especially when the signal is weak.
 
In \cite{zou2019quickestnetwork}, the problem of sequentially detecting a significant change (i.e., when at least $\eta$ number of sensors are affected by the change) was investigated. The event is dynamic, i.e., different nodes are affected at different times. Instead of using a scan statistic, which is computationally costly, a spartan-CUSUM (S-CUSUM) algorithm was constructed, which compares the sum of the smallest $N-\eta+1$ local CUSUM statistics to a threshold, where $N$ is the total number of nodes. For the case where the change propagates along network edges, a network-CUSUM (N-CUSUM) algorithm was further constructed based on the idea that the affected nodes shall induce a connected subgraph. The N-CUSUM algorithm was also shown to be first-order asymptotically optimal, and performs much better than the S-CUSUM numerically. The decentralized setting where there is no fusion center and nodes can only communicate with their neighbors was studied in \cite{zou2018icassp,li2019asilomar}, and the approach is based on a novel combination of the alternating direction method of multipliers (ADMM) and average consensus approaches. In \cite{kurt2018multisensor}, a Bayesian approach is used to model the dynamic change with an unknown propagation pattern, where the goal is to detect the change when it firstly emerges in the network; an optimal solution structure is derived using a dynamic programming framework.  

\subsection{Robust Sequential Change Detection}

Many classical procedures (for instance, CUSUM and SR) require exact knowledge of the pre- and post-change distributions. However, in real-world scenarios, the actual data distributions may be complex and different from what we have assumed. There can be adversarial attacks that significantly perturb the data distributions. This can lead to performance degradation of the optimal procedures. How to make the procedures more robust in the presence of model mismatch is the topic of robust sequential change detection. 

\subsubsection{Robustness to Model Uncertainties}\label{sec:robust}
There have been many efforts to make the detection procedure more robust to model uncertainties. One approach is to treat the pre- and post-change distributions to belong to some parametric family with unknown parameters in {\it uncertainty sets} and then form the GLR based test as we discussed earlier in Section~\ref{sec:other_seq_ch}. Another approach to developing good tests in the presence of model uncertainties is through the use of minimax robustness as the criterion as is done in the seminal work of Huber on robust hypothesis testing \cite{huber1965robust,huber2011robust}. The solution to the robust hypothesis testing problem usually relies on finding the least favorable distributions (LFDs) within the uncertainty classes, with likelihood ratio of these distributions used in constructing the robust tests. It can be shown that LFDs exist for uncertainty classes satisfying a certain joint stochastic boundedness (JSB) condition \cite{veeravalli1994minimax}. The problem of minimax robust sequential change detection was explored in  \cite{unnikrishnan2011minimax}, in which an exactly optimal solution was obtained for uncertainty classes satisfying the JSB condition under a generalized Lorden criterion. An extension of this result to asymptotic minimax robust sequential change detection is studied in \cite{molloy2017misspecified}, where a weaker notion of stochastic boundedness is introduced. 

A robust CUSUM algorithm is developed in \cite{cao2017robust}  by making a connection to convex optimization, which is particularly useful for the high-dimensional setting and leads to a tractable formulation. For instance, assuming the covariance matrix lies in an {\it uncertainty set} centered around a nominal value, the problem of finding LFDs can be cast as solving a semidefinite program and can be solved efficiently.  

\subsubsection{Robustness to Adversarial Attacks} 

The problem of sequential change detection in sensor networks in the presence of adversarial attacks \cite{lamport2019byzantine} was investigated in \cite{fellouris2017efficient,bayraktar2015byzantine}. In the presence of Byzantine attacks, an adversary may modify observations arbitrarily to defer the detection of a change and increase the false alarm rate. In \cite{bayraktar2015byzantine}, it is assumed that the change affects all but one compromised sensor, and the detection strategy is to raise a global alarm until two local CUSUMs exceed the threshold. In \cite{fellouris2017efficient}, a more general setting was investigated, where an unknown subset of sensors can be compromised. Sequential detection strategies were designed by waiting until $L$ local CUSUM statistics exceed the threshold (simultaneously or not) or by comparing the sum of the $L$ smallest CUSUM statistics to a threshold. With a proper choice of $L$, the above approaches are robust to Byzantine attacks. 

%

\subsection{Data-efficient Sequential Change Detection}

There is usually a cost associated with making observations in practical engineering applications, e.g., the power consumption in sensor networks. 
An extension of Shiryaev's formulation (Section~\ref{sec:bayes_optimaly}) was investigated in \cite{bane-veer-sqa-2012} by including an additional constraint on the average number of observations taken before the change. The cost of observations after the change is included in the detection delay. Specifically, whether to take an observation at time $t$ is controlled by an on-off binary control variable $S_t$, and $S_t$ is a function of all the information available up to time $t-1$. A data-efficient Shiryaev (DE-Shiryaev) algorithm was constructed in \cite{bane-veer-sqa-2012}, and was shown to be asymptotically optimal as $\PFA$ goes to zero.  The DE-Shiryaev  algorithm is also shown to have good observation cost-delay tradeoff curves: for moderate values of $\PFA$, for Gaussian observations, the delay of the algorithm is within 10\% of the Shiryaev delay even when the observation cost is reduced by more than 50\%.
Furthermore, the DE-Shiryaev algorithm is substantially better than the standard approach of \textit{fractional sampling} scheme, where the Shiryaev algorithm is used and where the observations to be skipped are determined \emph{a priori} in order to meet the observation constraint. A minimax formulation was further proposed in \cite{banerjee2013data} to address the scenario when a prior on the change-point is not available. 
The DE-CUSUM algorithm developed in \cite{banerjee2013data} is shown to be asymptotically optimal as $\FAR$ goes to zero,  and significantly outperforms fractional sampling in simulations. 
Extensions to composite post-change distributions were studied in \cite{banerjee2015data}, and generalizations to distributed sensor networks were explored in \cite{banerjee2015data2}.

\subsection{High-dimensional Streaming Data}

High-dimensional data usually have low-dimensional structures, such as sparsity and low-rankness, which can be leveraged to achieve improved detection performance and computational efficiency. Meanwhile, missing data is very common for high-dimensional streaming data. In this section, we review recent advances in these directions.

\subsubsection{Sparse Change in Multiple Data Streams}\label{sec:sparse-multi-stream}

For multiple independent streams of data, a mixture procedure was developed in \cite{xie2013sequential} to monitor parallel streams for a  change-point that affects only a subset of them (usually sparse). Both the subset being affected and the post-change distribution are unknown. The mixture model hypothesizes that each sensor is affected with a small probability $\varrho \in (0, 1)$ by the change, where $\varrho$ is pre-specified. The mixture detection statistic at time $t$ is defined as \[\sum_{n=1}^N \log\left[1-\varrho + \varrho f_1(X_{t}^{(n)})/f_0(X_{t}^{(n)})\right],\]
where $X_{t}^{(n)}$ denotes the observation at the $n$-th sensor and at time $t$, and $N$ is the number of sensors. Another efficient global monitoring scheme was proposed in \cite{mei-ieeetit2015} by combining hard thresholding with linear shrinkage estimator for the post-change parameters. In recent works \cite{zhang2018asymptotic,liu2019scalable}, a similar problem was tackled by running local detection procedures and using the sum of the shrinkage transformation of local detection statistics as a global detection statistic. 
This sum-shrinkage framework was further extended in \cite{zhang2019robust} to be more robust to outliers using the Box-Cox transformation. Recent work \cite{enikeeva2019high} studied change detection in regimes where the dimension tends to infinity and the length of the sequence grows with the dimension. 

\subsubsection{Subspace Change Detection}

In many applications, the change in high-dimensional data covariance structure can be represented as a low-rank change. For instance, in seismic signal detection \cite{xie2020seq_ana}, a similar waveform is observed at a subset of sensors after the change. Such a change can be modeled as the covariance matrix shifts from an identity matrix to a ``spiked'' covariance model \cite{johnstone2001distribution}. The subspace-CUSUM procedure was developed in \cite{xie2020seq_ana}, in which the {\it unknown} subspace in the post-change spiked model is estimated sequentially and further used to obtain the log-likelihood ratio statistic. A CUSUM procedure for detecting switching subspace (from a known subspace to another target subspace) was studied in  \cite{Gu2018Subspace}.

\subsubsection{Missing Data}

In high-dimensional time series, it is common that we cannot observe all the entries at each time. The missing components in the observed data handicap conventional approaches. In \cite{mousse}, a mixture type of approach was proposed by combining subspace tracking with missing data to model the underlying dynamic of data geometry (submanifold). Specifically, streaming data is used to track a submanifold approximation, to measure deviations from this approximation, and to calculate a series of statistics of the deviations for detecting when the underlying manifold has changed.

\subsubsection{Sketching to Conquer High-dimensionality}

To detect changes quickly over high-dimensional data, we may need to conquer the challenges presented by the data's high dimensionality. {\it Sketching} is a commonly used strategy to reduce data dimensionality, which performs linear projections of high-dimensional data into a small number of sketches. A GLR procedure based on data sketches was studied in \cite{sketchingChange2015}, with the precise characterization of performance metrics and the minimum number of sketches needed to achieve good performance. Multiple types of sketching matrices can be used, such as Gaussian random matrices, expander graphs, and network topology constrained matrices. The sketching procedure is relevant to large power networks where we cannot place a sensor on each node or edge. Instead, each sensor will measure {\it aggregates} of the network states at a few edges or nodes. In \cite{sketchingChange2015}, the mean-shift detection problem in power networks is studied, where each measurement corresponds to a linear combination of the state at an edge, e.g., real power flow. This leads to a sketching matrix determined by the network topology. 

\subsection{Joint Detection and Estimation}

It is common that the distribution after the change is unknown. For instance, before the change in industrial process monitoring applications, the production line is in-control and well-calibrated (thus the distribution before the change is known). However, after the change, an anomaly causes a shift to the operation into an unknown status. Therefore, it is interesting to incorporate estimates of the possible post-change status into the detection statistic when performing detection; this problem is related to robust sequential change detection, as discussed in Section\,\ref{sec:robust}. In other situations, we need to estimate the post-change distribution in retrospect for identifying the change. There has been much work establishing the theoretical foundation for joint detection and estimation. For instance, \cite{moustakides2012joint} combines the Bayesian formulation of the estimation and detection and develops an optimal procedure to achieve a tradeoff between detection power and estimation quality.
In another context, it is also referred to as sequential change diagnosis \cite{dayanik2008bayesian}. Quickest {\it searching} of the change-point (e.g., quickest search for rare events) has been developed in \cite{tajer2013quick,heydari2018quickest,heydari2019quickest}.

\subsection{Spatio-temporal Change Detection}

When modeling discrete event data, the point process model \cite{daley2007introduction} is frequently used due to its capability of modeling the time intervals between events directly. Point processes assume that time intervals between events are exponentially distributed. For example, in Poisson processes the intervals are independent, and in Hawkes processes the intervals are dependent, and the intensity depends on the events that occurred in the past \cite{coevolve2015}. The ``autoregressive'' nature of Hawkes processes makes them attractive in modeling temporal dependence and causal relationships, including market models \cite{Toke2010}, earthquake event prediction \cite{Ogata88}, inferring leadership in e-mail networks 
\cite{fox2014modeling}, and topic models \cite{HawkesTopic15}. The multi-dimensional Hawkes process model over networks can model highly correlated discrete event data \cite{reinhart2018review} and capture dependence over networks and propagation of the signal in such settings. 

Detection of changes for point processes has attracted much attention for both single event stream and multiple streams over networks (or over multiple locations). For example, there are works focusing on Poisson processes \cite{ShenZhang12,zhang2014scanning,herberts2004optimal}, and some recent work on one-dimensional \cite{Ludkovski12,pinto2015trend} and multi-dimensional (network) point processes \cite{li2017detecting,wang2020detecting}. In particular, \cite{li2017detecting} studied the change detection for networked streaming event data and constructed GLR type procedures; \cite{wang2020detecting} developed the penalized dynamic programming algorithm to detect coefficient changes in discrete-time high-dimensional self-exciting Poisson processes in an offline setting.

This topic is also related to the multisource quickest detection problem, mostly assuming independence between multiple data streams. For instance, the quickest detection of the minimum of change-points for two independent compound Poisson processes was considered in  \cite{bayraktar2007quickest} and optimal Bayesian sequential detection procedures were developed. 

\subsection{Change Detection-Isolation-Identification}
In addition to detecting the change quickly after it occurs, sometimes we are also interested in identifying the post-change model and/or isolating a subset of nodes within a large network affected by the change. In \cite{dayanik2013asymptotically}, an asymptotically optimal Bayesian detection–isolation scheme was proposed assuming the post-change model is one of the finitely many distinct alternatives. In a series of works, Nikiforov introduced a minimax optimal detection-isolation algorithm for stochastic dynamical systems \cite{nikiforov1995generalized}, developed a recursive variant of the algorithm that achieves better computational efficiency \cite{nikiforov2000simple}, and provided an asymptotic lower bound for the mean detection-isolation delay with constraints on the probability of false isolation and the average time before a false alarm \cite{nikiforov2003lower}. Natural generalizations of CUSUM and SR procedures for detection-isolation problems were discussed in \cite{tartakovsky2008multidecision}. See  \cite{tartakovsky2014sequential,tartakovsky2019sequential} for more detailed overviews.

\subsection{Alternative Performance Metrics}

Other than what have been presented in this survey, many alternative performance metrics have also been considered. For instance, \cite{poor1998quickest} investigated an exponential penalty of delay rather than a linear penalty (as used in the definition of $\CADD$, for instance). Such performance measures can be more accurate, sometimes for financial applications. In these cases, the change-point may not represent a time at which a fundamental shift in the performance occurs, but the compounding of investment growth can be a more suitable measure of the cost of delay. Similarly, in the health monitoring of components in aircraft systems, communication networks, and power grids,  the effects of undetected faults can exponentiate with time. For problems involving estimation, the performance measures can also involve estimation accuracy, for instance, change-point location and other parameters involved in the problem. With many parallel data streams, the error metric can be the false discovery rate ($\FDR$), which is the expected ratio of the number of falsely declared data streams to the total number of declared data streams \cite{chen2020false}.

\section{New dimensions}\label{sec:new-dimension}

\subsection{Machine Learning and Change Detection} 

Modern machine learning approaches can be adopted for solving sequential change detection problems, which we will review in this subsection. 

\subsubsection{Density Ratio Estimation} 

Instead of estimating the post-change density $f_1$ as in the GLR procedure, we may estimate the density ratio $f_1/f_0$ directly (referred to as density ratio estimation \cite{sugiyama2012density}), based on which we develop sequential change detection procedures. A data-driven framework using neural networks was developed in \cite{moustakides2019training}. More specifically, given two sets of data sampled from the densities of interest, an optimization problem is defined so that the solution, specified through neural networks, will correspond to the desired likelihood ratio function or its transformations and can then be used for sequential change detection.  

\subsubsection{Anomaly Detection} 

Change detection is closely related to {\it anomaly detection}, which is a popular topic in machine learning and data mining, and many machine learning techniques have been developed. In particular, an recurrent neural network (RNN) based approach computes the detection statistic (referred to as the anomaly score) in an online fashion and compares with a threshold for anomaly detection \cite{Saurav2018ACM}. The RNN-based approach can benefit certain situations since they are known to capture complex temporal dependencies for multivariate time series. We refer to \cite{chalapathy2019deep} for a recent survey on deep learning techniques for anomaly detection. Developing mathematical theory for RNN-based sequential change detection is still an open question. 

\subsubsection{Online Learning and Change Detection}

Online implementation is one of the most critical aspects of sequential change detection algorithms in practice. Although many algorithms enjoy recursive structure, such as CUSUM and SR procedures, some sequential detection procedures face a significant hurdle of online implementation due to their non-recursive nature. For instance, window-limited GLR statistic, although enjoying robust performance in the presence of unknown post-change distributions, is not recursive since the parameters need to be continuously estimated by incorporating new samples.  To tackle this challenge, inspired by online learning, \cite{cao2018entropy} develops an online mirror descent-based GLR procedure to update the estimate of the unknown post-change parameter with new data. Another highly cited work \cite{adams2007bayesian} develops an online change detection procedure based on Bayesian computing. In recent work, \cite{titsias2020sequential} develops a framework for joint sequential change detection and online model fitting, which will be particularly suitable for parameterized models. A GLR procedure is developed in this framework using estimates of the unknown high-dimensional parameter obtained by the gradient descent update. 

\subsubsection{Tracking Data Dynamics} 

Many sequential data are dynamic even before the change has happened; for instance, solar flare detection from satellite video streaming \cite{Xie2012,mousse}. To build methods that work with real-world scenarios, we need to develop robust methods that can adapt to normal data dynamics without mislabeling them as change-points. A possible strategy is to combine tracking with detection. For instance,  \cite{Xie2012,mousse} developed a procedure to detect sparse changes when the pre-change high-dimensional data is time-varying. The data dynamic is captured by tracking a time-varying manifold using variants of subspace tracking (e.g., GROUSE \cite{GrouseGLOBALconv2015}, PETRELS \cite{ChiPETRELS13}, or MOUSSE algorithm \cite{mousse}). Another instance is the network Hawkes process model, where we may track the Hawkes process through online learning techniques \cite{HallWillett2014}.

\subsubsection{Active Learning and Change Detection}

For certain applications such as material science and recovering seafloor depth, data acquisition is expensive. Thus, it is desirable to collect data that is most useful in a sequential fashion, which is the theme of active learning (see, e.g., \cite{singh2006active,castro2008active}). The combination of active learning and change detection was introduced as active change-point detection (ACPD) problem in \cite{hayashi2019active}. The task is to adaptively determine the next input to detect the change-point in a black-box expensive-to-evaluate function, with as few evaluations as possible. The method utilizes the existing change detection method to compute change scores and a Bayesian optimization method to determine the next input. A CUSUM procedure with an adaptive sampling strategy to detect mean shifts was developed in  \cite{liu2015adaptive}.
%

\subsubsection{Detection with Data Privacy}  

As data privacy has growing importance in modern applications in social settings, it also leads to developing private change detection algorithms. Both offline and online change detection methods through the lens of {\it differential privacy} have been developed in \cite{cummings2018differentially}. A different privacy-aware sequential change detection method was studied in \cite{lau2020privacy},  using {\it maximal leakage} as the privacy metric, which is a weaker form of privacy compared with \cite{cummings2018differentially}.

\subsubsection{Change Detection for Reinforcement Learning}  

Reinforcement learning is a major type of sequential decision-making methodology in the era of artificial intelligence. How to implement reinforcement learning in a non-stationary and changing environment is still a mostly unexplored area. Recently, there have been some attempts to combine sequential change detection and reinforcement learning \cite{padakandla2019reinforcement}, where change detection algorithms are utilized to detect the transition of the environment and trigger transitions of reinforcement learning algorithms.  

\subsection{Distribution-free Methods}

Distribution-free methods aim to detect the change without making explicit distributional assumptions on the data. Such methods are particularly attractive in machine learning, such as kernel MMD based method discussed in Section\,\ref{change-of-measure}, due to their flexibility in working with complex data. There have been kernel-based non-parametric methods developed in terms of change detection, both for the offline setting \cite{harchaoui2009kernel,harchaoui2013kernel,arlot2019kernel} and the online setting \cite{li2019scan}. MMD statistics have also been used for anomalous sequence detection, for instance, \cite{7997789,zou2017nonparametric}. Besides MMD, other distribution-free methods have been developed for change detection. For instance, dissimilarity measures based on the kernel support vector machine (SVM) were built in \cite{desobry2005online}, and  generalized likelihood test directly using data empirical distributions when the true distributions are supported on a finite alphabet were constructed in \cite{nitinawarat2017universal,nitinawarat2015universal,lau2018binning,bu2019linear}.

There are many other types of distribution-free non-parametric tests for change detection developed in various contents. For instance, the maximal $k$-largest sample coherence between columns of each observed random matrix was developed to detect change for large-scale random matrices \cite{banerjee2018quickest}. A nearest-neighbors-based statistic was proposed in \cite{chen2019sequential} to detect the change in sequences of multivariate observations or non-Euclidean data objects such as network data. The weighted moving averages were studied in \cite{frias2014online} to detect univariate drifts. 
A non-parametric approach was developed in \cite{pawlak2013nonparametric} to detect departure from the reference signal with non-i.i.d. underlying time series. The spectral scan statistic for change detection over graphs was considered in \cite{pmlr-v31-sharpnack13b}.  Wasserstein distance was used to detect segments of times series in \cite{cheng2020optimal}. In \cite{balsubramani2015sequential}, test statistics were constructed using martingales under the null hypothesis, and the rejection threshold is determined using a uniform non-asymptotic law of the iterated logarithm. 

\subsection{Non-stationary Multi-armed Bandits with Changes}

Multi-armed bandit is a class of fundamental problems in online learning and sequential decision-making. A learning agent aims to maximize its expected cumulative reward by repeatedly selecting to pull one arm at each time step. Change detection can play a role in the scenario where the reward distributions may change in a piece-wise-stationary fashion at unknown time steps. To handle dynamic multi-armed bandit problems, various change detection methods were considered, including the Page-Hinkley test  \cite{hartland:inria-00164033}, a windowed mean-shift detection \cite{Yu-Mannor-2009-ICML}, CUSUM test \cite{AAAI18BanditCD}, and sample mean based test \cite{cao2019bandit}. Usually,  the algorithm will reset once a change is detected. From a Bayesian perspective, the Thompson sampling strategy equipped with a Bayesian change-point mechanism was considered in \cite{pmlr-v31-mellor13a}. The adversarial multi-armed bandit problem with change points was also considered in \cite{DSAA-2015-Bandits}.

\subsection{Optimization for Change Detection and Estimation}

Optimization is becoming a centerpiece in developing modern machine learning algorithms. Recent advances in convex optimization have enabled solving many large-scale problems. A line of research aims to casts (offline) change detection and estimation (of their locations) as an optimization problem. The benefits of this optimization-based approach typically include computational efficiency (when the optimization problem is convex) and theoretical performance guarantees based on optimization theory. Below we give some examples.

The univariate change detection for a mean shift using an optimization approach has been studied in \cite{lin2016approximate}, and performance guarantees were established by relating the $\ell_2$ recovery error to detection performance. 
A $\ell_0$-penalized least squares method was considered in \cite{wang2020univariate}. By connecting to binary segmentation methods, change detection and localization for univariate data in the non-parametric settings was studied in \cite{padilla2019optimal}.

Multivariate change detection using an optimization approach has also been studied. For instance, a dynamic programming approach was developed for recovering an unknown number of change-points from multivariate autoregressive models \cite{wang2019localizing}. A network binary segmentation method for change detection was proposed in \cite{wang2018optimal}, which has been extended for covariance matrix change detection in \cite{wang2017optimal}. Finally, the work \cite{soh2017high} combined the filtered derivative with convex optimization methods to estimate change-points for multi-dimensional data.

\section{Modern applications}\label{sec:applications}

Sequential change detection has traditionally been used in industrial process monitoring applications, which was probably the original motivation for change detection procedures to be developed in the early days. The wide adoption of change detection in industrial quality engineering and manufacturing initiates the field of statistical process control (SPC) (see, e.g., \cite{shi2006stream, oakland2018statistical}). Recently, there have been many more modern applications for sequential change detection, and we present a selection of them here. 



\subsection{Smart Grids}

The sequential change detection methodology has been recently successfully applied for sequential line outage detection in power transmission systems. In modern smart grids, high-speed synchronized voltage phase angle measurements are taken from phasor measurement units (PMU). Based on PMU measurements, a linearized incremental small-signal power system model was developed in \cite{chen2015quickest}. Once a line outage occurs, there is a change in the covariance matrix of incremental phases, by monitoring which, line outages can be detected and identified using sequential change detection algorithms. In \cite{Rovatsos2:2016}, the transient dynamics of the power system following a line outage is further incorporated. The D-CUSUM algorithm was then developed to incorporate the dynamic nature of the line outage in \cite{Rovatsos2:2016} (see Section \ref{sec:transient} for more details).

There have been other works on sequential change detection for smart grids. The generalized local likelihood ratio test was applied for voltage quality monitoring \cite{li2015cooperative}, photovoltaic  systems \cite{chen2016quickest}, attack detection in the multi-agent reputation systems \cite{li2014quickest}, wide-area monitoring \cite{Li2015smart}, and cyber-attacks detection in discrete-time linear dynamic system \cite{Kurt2018grids,Kurt2019grids}. The decentralized detection with level-triggered sampling was considered in  \cite{Yilmaz2013Channel}. In \cite{heydari2017quickest}, a general stochastic graphical framework for modeling the bus measurements and a data-adaptive data-acquisition and decision-making processes were designed for the quickest search and localization of anomaly in power grids.

\subsection{Cybersecurity}
Cybersecurity has become a critical problem with the development of wireless communication, networking, and the Internet of things. It is of practical importance to detect attacks and intrusions in real-time from network streaming data, e.g., denial-of-service attacks, worm-based attacks, port-scanning, and man-in-the-middle attacks. The sequential change detection approach is a natural fit since the attacks usually change network traffic distribution. In \cite{tartakovsky2006detection}, multi-channel generalizations of the CUSUM procedure and non-parametric tests were proposed. In \cite{tartakovsky2006novel}, adaptive sequential methods were proposed for early detection of subtle network attacks, utilizing data from multiple layers of the network protocol. In \cite{tartakovsky2012efficient}, a multi-cyclic detection procedure based on the SR procedure was proposed. In \cite{tartakovsky2014rapid}, score-based CUSUM and SR procedures were exploited for network anomaly detection, and a hybrid detection system 
was proposed. The application to cybersecurity was also discussed in books \cite{tartakovsky2014sequential,tartakovsky2019sequential}, and recent reviews \cite{jeske2018statistical,jeske2014statistical}.

\subsection{Sensors Networks} 
Sensor networks collecting sequential data have been widely used for geophysical, environmental, traffic, and internet traffic monitoring applications, which we will briefly summarize in this subsection. 

Seismology is experiencing rapid growth in the quantity of data. Earthquake detection aims to identify seismic events in continuous data -- a fundamental operation for seismology \cite{yoon2015earthquake}. Modern ultra-dense seismic sensor arrays have obtained a massive amount of continuous data for seismic studies, and many such data are publicly available through IRIS \cite{IRIS}. In the old days, network seismology treated seismic signals individually - one sensor at a time - and detected an earthquake upon multiple impulsive arrivals consistent with a source within the Earth \cite{Johnson94}. Recently, with advances in sensor technology, which bring densely sampled data, high-performance computing and high-speed communication, we are able to use a {\it network-based detection} by exploiting correlations between sensors to extract coherence signals. This will enhance the systematic detection of weak and unusual events that currently go undetected using individual sensors. Detecting such weak events is very crucial for earthquake prediction \cite{earthquake_event_70,forecasting_weak_events2013}, oil field exploration, volcano monitoring, and deeper earth studies \cite{FanchiScience15}. Towards this goal, in \cite{xie2020seq_ana}, a subspace-CUSUM procedure was developed for network-based detection by exploiting the low-rank subspace structure induced by waveform similarity.

Sensor networks have also been deployed to monitor drinking water safety from the water tower to private residences. Sequential change detection using residual chlorine concentration measurements from the sensors network was developed in \cite{guepie2012sequential}. Methods have also been developed for monitoring river contamination \cite{chen2019reduce,chen2017textsf}, which specifically consider the spatio-temporal correlation in observations along the sensor network due to water dynamics.  

Sequential monitoring of traffic flow using traffic sensors has been considered in \cite{rajagopal2008distributed}, and a distributed, online, sequential algorithm for detecting multiple faults in a sensor network was presented therein. Recently,   Hawkes processes models for correlated traffic anomalies using data collected by inductive-loop traffic detectors were developed in \cite{zhu2020traffic}. 


\subsection{Wireless Communications}

Sequential change detection has been used for wireless communications, including online user activity detection for multi-user direct-sequence/code-division multiple-access (DS-CDMA) environment   \cite{oskiper2002online}, detecting ``spectrum opportunities'' in the cognitive radio setting by identifying the occupancy and idle of channels from primary user's activities \cite{LaiCognitive08,xie2012spectrum}. More recently, \cite{huang2020lpd} established a change detection framework for low probability of detection (LPD) communication, where a transmitter, Alice, wants to hide her transmission to a receiver, Bob, from an adversary, Willie; three different sequential tests were considered, including  Shewhart, CUSUM, and SR procedures, to model Willie's detection process.

\subsection{Video Processing and Computer Vision}

\begin{figure}[!ht]
\centering
\includegraphics[width=0.9\linewidth]{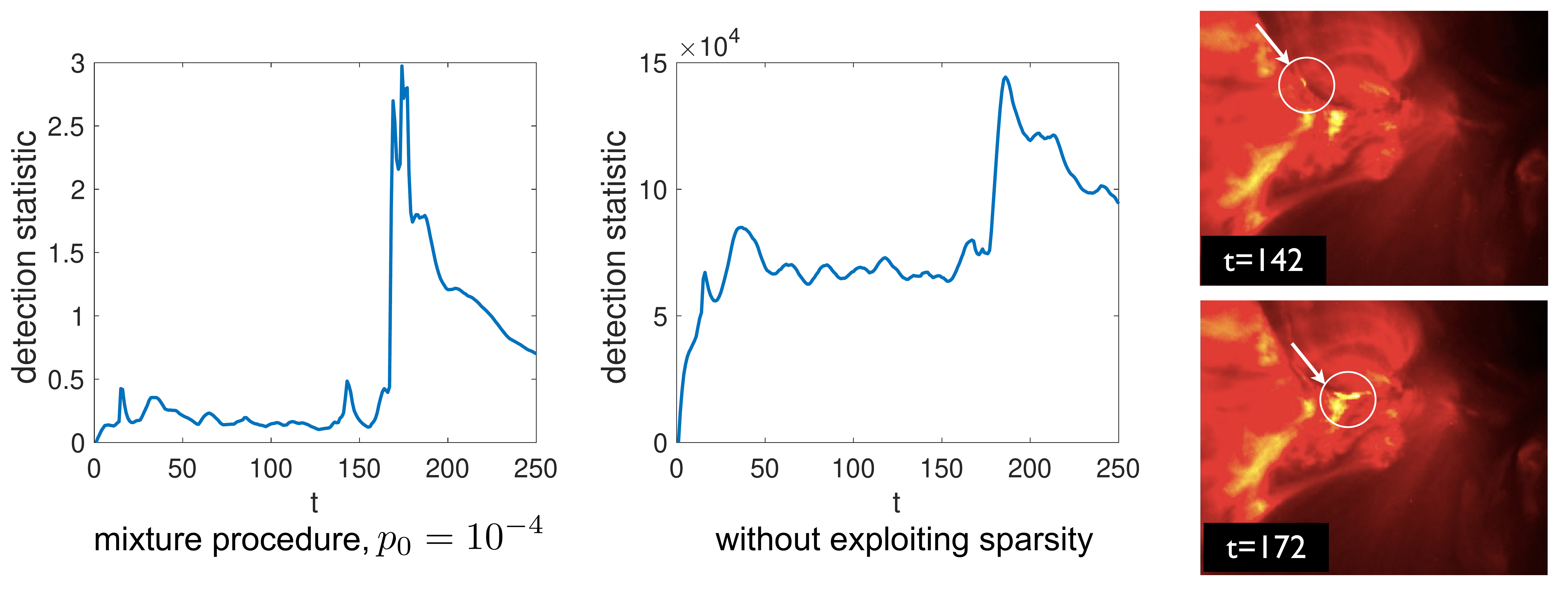}
\caption{Solar flare detection with the mixture procedure as considered in \cite{mousse}; the first minor solar flare at $t = 142$ is hardly visible, and it is missed entirely by a baseline detection statistic (sum of CUSUM at each data dimension without exploiting sparsity in the change: ``solar flare''). This also illustrates the importance of exploiting sparsity in the change.}
\label{fig:solar}
\end{figure}

Change detection is one of the most commonly encountered low-level tasks in computer vision and video processing \cite{radke2005image}, and many such problems are essentially sequential. A plethora of practical algorithms have been developed to date; for instance, scene change detection \cite{lelescu2003statistical}, street-view change detection \cite{alcantarilla2018street}, and change detection in video sequences
\cite{tsechpenakis2006robust}. In \cite{de2014change}, a pixel-based weightless neural network (WNN) method was developed to detect changes in the field of view of a camera. In \cite{lingg2014sequential}, multiple images from reference and mission passes of a scene of interest were used to improve  detection performance. There are still many open questions regarding how to leverage the power of statistical sequential change detection for computer vision and video processing. We present an example of solar flare detection from video sequences in Fig.\,\ref{fig:solar}, which has been considered in several works along this line including \cite{mousse}. 


\subsection{Social Networks}

The wide-spread use of social networks and the great availability of information networks (e.g., Twitter, Facebook, blogs) lead to a large amount of user-generated data \cite{krishnamurthy2014tutorial}, which are quite valuable in studying many social phenomena. One important aspect is to detect change-points in streaming social network data \cite{krishnamurthy2012quickest}, which may represent the collective anticipation of or response to external events or system ``shocks'' \cite{PeelClauset2014}. Detecting such changes could provide a better understanding of the patterns of social life. In other cases, early detection of change-points can predict or even prevent social stress due to disease or international threat, for instance, detecting self-exciting changes (modeled by network Hawkes processes) in social networks \cite{li2017detecting}. A related topic is distributed hypothesis testing in social networks: \cite{lalitha2018social} showed the exponential convergence rate of a Bayesian update scheme of nodal belief (distribution estimate) in the social learning setting.   
%


\subsection{Epidemiology}

Sequential change detection can potentially play an important role in public health and disease surveillance. Early detection of epidemics is a very important topic. 
In \cite{baron2002bayes,epidemic_change_04}, Baron cast the early detection of epidemics as a Bayes sequential change detection problem and proposed an asymptotically pointwise optimal stopping rule, which is computationally efficient for complicated prior distributions arising in epidemiology. In \cite{yu2013change}, a modified CUSUM procedure was proposed for the susceptible–infected–recovered (SIR) epidemic model to detect change-point in the infection rate parameter. Moreover, change detection has been incorporated into studying the intervention's effectiveness, based on the premise that the underlying epidemiological model may change over time due to interventions. Evaluating intervention measures' effectiveness requires detecting underlying change-points, which becomes even more important in the COVID-19 era \cite{zhu2020high}. Such works include  \cite{wieland2020phenomenological, pedersen2020data}, which estimate the change-points in time series to assess the effectiveness of interventions such as lock-down and mask usage;  in \cite{dehning2020inferring},  the problem of detecting the growth rate change for the COVID-19 spread in Germany was studied, where results were further incorporated into forecasting. There are still many open questions in this area regarding developing effective sequential change detection procedures suitable for infectious disease early detection.



\section{Conclusions} \label{sec:discussion}

Our goal in this survey was to provide a glimpse of the past and recent advances in sequential change detection,  and its application in various domains. We have covered different types of sequential change detection procedures, both theoretically optimal and practical. We also discussed how the intersection of sequential change detection with other areas has created interesting new directions for research. 

\section*{Acknowledgement}
The authors are grateful to the Guest Editor and the anonymous reviewers for their helpful comments.

\bibliographystyle{IEEEtranS}
\bibliography{full_ref}

\end{document}